\documentclass[draft]{amsart}
\let\oldtocsection=\tocsection

\let\oldtocsubsection=\tocsubsection

\let\oldtocsubsubsection=\tocsubsubsection

\renewcommand{\tocsection}[2]{\hspace{0em}\oldtocsection{#1}{#2}}
\renewcommand{\tocsubsection}[2]{\hspace{1em}\oldtocsubsection{#1}{#2}}
\renewcommand{\tocsubsubsection}[2]{\hspace{2em}\oldtocsubsubsection{#1}{#2}}

\usepackage[Symbol]{upgreek}

\usepackage{amssymb}
\usepackage{epic}
\usepackage{mathrsfs}
\usepackage{accents}
\usepackage{stmaryrd}
\usepackage{enumitem}
\usepackage{amsmath}
\usepackage{graphicx}
\usepackage{lipsum}
\usepackage{braket}
\usepackage{hyperref}
\usepackage{mathtools}
\usepackage{dsfont}
\usepackage{enumitem}   
\DeclarePairedDelimiter{\ceil}{\lceil}{\rceil}
\DeclarePairedDelimiter{\floor}{\lfloor}{\rfloor}
\usepackage{kantlipsum} 

\calclayout

\input{xy}
\xyoption{matrix}
\xyoption{arrow}

\numberwithin{equation}{section}

\renewcommand\epsilon{\varepsilon}

\def \<{\langle}
\def \>{\rangle}

\def \tilde {\widetilde}

\def \hat {\widehat}

\def \((  {(\!(}
\def \)) {)\!)}

\DeclareMathSymbol{\precequ}{\mathrel}{symbols}{"16}
\DeclareMathSymbol{\succequ}{\mathrel}{symbols}{"17}

\newtheorem{theorem}{Theorem}[section]
\newtheorem{lemma}[theorem]{Lemma}
\newtheorem{prop}[theorem]{Proposition}
\newtheorem{cor}[theorem]{Corollary}

\newtheorem*{theorem*}{Theorem}

\theoremstyle{definition}

\theoremstyle{remark}

\newtheorem*{remark}{Remark}

\newtheorem*{conjecture}{Conjecture}

\def \sgn {\operatorname{sign}}

\def \GL{\operatorname{GL}}

\let\oldi\i
\let\oldj\j
\renewcommand\i{\relax\ifmmode{\boldsymbol{i}}\else\oldi\fi}
\renewcommand\j{\relax\ifmmode{\boldsymbol{j}}\else\oldj\fi}

\renewcommand\leq{\leqslant}
\renewcommand\geq{\geqslant}

\DeclareMathAlphabet{\mathbf}{OML}{cmm}{b}{it}

\DeclareFontFamily{U}{fsy}{}
\DeclareFontShape{U}{fsy}{m}{n}{<->s*[.9]psyr}{}
\DeclareSymbolFont{der@m}{U}{fsy}{m}{n}
\DeclareMathSymbol{\der}{\mathord}{der@m}{182}

\DeclareSymbolFont{der@m}{U}{fsy}{m}{n}
\DeclareMathSymbol{\derdelta}{\mathord}{der@m}{100}

\DeclareSymbolFont{imag@m}{OT1}{cmr}{m}{ui}
\DeclareMathSymbol{\imag}{\mathord}{imag@m}{105}

\DeclareFontFamily{OMS}{smallo}{}
\DeclareFontShape{OMS}{smallo}{m}{n}{<->s*[.65]cmsy10}{}
\DeclareSymbolFont{smallo@m}{OMS}{smallo}{m}{n}
\DeclareMathSymbol{\smallo}{\mathord}{smallo@m}{79}

\DeclareFontFamily{OMS}{largerdot}{}
\DeclareFontShape{OMS}{largerdot}{m}{n}{<->s*[.8]cmsy10}{}
\DeclareSymbolFont{largerdot@m}{OMS}{largerdot}{m}{n}
\DeclareMathSymbol{\largerdot}{\mathord}{largerdot@m}{15}

\newcommand\Aut{\operatorname{Aut}}

\newcommand\diag{\operatorname{diag}}

\DeclareMathSymbol{\llambda}{\mathord}{der@m}{108}
\DeclareMathSymbol{\rrho}{\mathord}{der@m}{114}

\newcommand{\equationqed}[1]{\[\pushQED{\qed}#1 \qedhere\popQED\]\let\qed\relax}
\newcommand{\alignqed}[1]{\begin{align*}\pushQED{\qed} #1 \qedhere\popQED\end{align*}\let\qed\relax}

\makeatletter
\newcommand{\dminus}{\mathbin{\text{\@dminus}}}

\newcommand{\@dminus}{%
  \ooalign{\hidewidth\raise1ex\hbox{\bf.}\hidewidth\cr$\m@th-$\cr}%
}
\makeatother

\DeclareMathOperator{\SL}{SL}

\DeclareMathOperator{\J}{J}
\DeclareMathOperator{\SO}{SO}

\DeclareMathOperator{\PGL}{PGL}
\DeclareMathOperator{\Wd}{Wd}
\DeclareMathOperator{\JL}{JL}
\DeclareMathOperator{\Sym}{Sym}
\DeclareMathOperator{\Ad}{Ad}
\DeclareMathOperator{\Lie}{Lie}
\DeclareMathOperator{\Vol}{Vol}
\DeclareMathOperator{\fin}{fin}
\DeclareMathOperator{\hol}{hol}
\DeclareMathOperator{\He}{H}
\DeclareMathOperator{\X}{X}
\DeclareMathOperator{\Y}{Y}
\DeclareMathOperator{\Z}{Z}
\DeclareMathOperator{\N}{N}
\DeclareMathOperator{\Pe}{P}

\title{CENTRAL VALUES OF DEGREE SIX L-FUNCTIONS: THE CASE OF HILBERT MODULAR FORMS}
\author{Utkarsh Agrawal}

\begin{document}
\maketitle

\begin{abstract}
 In this article we will give a formula for the central value of the completed $L$-function $L(s,\Sym^{2} g\times f)$, where $f$ and $g$ are Hilbert newforms, by explicitly computing the local integrals appearing in the refined Gan-Gross-Prasad formula for $\SL_{2}\times\widetilde{\SL_{2}}$. We also work out the rationality of this value in some special cases and give a conjecture for the general case.    
\end{abstract}
\tableofcontents

\begin{section}{Introduction}\label{classlang}
Let $F$ be a totally real number field of degree $n$ over $\mathbb{Q}$. We denote by $\mathcal{O}, \mathbb{A},\Sigma_{\infty}$ respectively, the ring of integers of $F$, the ring of adeles of $F$ and the set of archimedean places of $F$. For an integral ideal $\mathfrak{n}\in\mathcal{O}$, let $\Gamma_{0}(\mathfrak{n})$ denote the congruence subgroup of $\GL_{2}^{+}(\mathcal{O})$, the group of $2\times 2$ matrices with entries in $\mathcal{O}$ and totally positive determinants, defined as
\begin{equation*}
    \Gamma_{0}(\mathfrak{n})=\Bigg\{\begin{pmatrix}
    a & b\\
    c & d 
    \end{pmatrix}\in \GL_{2}^{+}(\mathcal{O})
    :  
    c\in \mathfrak{n}\Bigg\}.
\end{equation*}
Let $f$ resp. $g$ be a holomorphic Hilbert newform of weight $2\kappa=(2\kappa_{v})_{v}$ resp. $\kappa'+1=(\kappa'_{v}+1)_{v}\in \mathbb{Z}^{\Sigma_{\infty}}$ (identified with $\mathbb{Z}^{n}$). We will assume that both $f$ and $g$ are of full level, that is, of level $\Gamma_{0}(\mathcal{O})$, and that $\kappa_{v}$ and $\kappa'_{v}$ are greater than 1 for all $v$. Let $\Sigma_{ub}$ and $\Sigma_{b}$ be subsets of $\Sigma_{\infty}$ whose union is $\Sigma_{\infty}$, such that $\kappa_{v}>\kappa'_{v}$ for $v\in \Sigma_{ub}$ and $\kappa_{v}\leq \kappa'_{v}$ for $v\in \Sigma_{b}$. Put
\begin{equation*}
    r_{v}=\begin{cases}
    \kappa_{v}-\kappa'_{v}-1 & \text{if}\ v\in \Sigma_{ub}\\
    \kappa'_{v}-\kappa_{v} & \text{if}\ v\in \Sigma_{b}.
    \end{cases}
\end{equation*}

Since the integral weight modular form $f$ is a Hecke eigenform, there is associated to it a unique (up to multiplication by a scalar) nonzero half-integral weight modular form, $h$, by the Shimura-Shintani-Waldspurger correspondence(cf. \cite{HI13}, \S9). The modular form $h$ is a holomorphic Hilbert newform and lies in the Kohnen plus space of modular forms of weight $\kappa+\frac{1}{2}$. Furthermore, this Kohnen plus space is isomorphic to a space of Jacobi forms (cf. \cite{HI13}, \S14). Thus to $h$ one can associate a Jacobi form, $F_{h}$, of weight $\kappa+1$ (and index 1), which, in particular, is a  holomorphic function on $\mathfrak{H}^{n}\times \mathbb{C}^{n}$. Write the coordinates on $\mathfrak{h}^{n}$ as $(\tau_{v})_{v}$ and put $$g_{ub}((\tau_{v})_{v})=g((\tau_{v})_{v\in\Sigma_{b}},(-\overline{\tau_{v}})_{v\in\Sigma_{ub}}),\ \ \  h_{ub}((\tau_{v})_{v})=h((\tau_{v})_{v\in\Sigma_{b}},(-\overline{\tau_{v}})_{v\in\Sigma_{ub}}).$$
Evidently, $g_{ub}$ (resp. $h_{ub}$) is the modular form on $\mathfrak{h}^{n}$ which is similar to $g$ (resp. $h$) but is \textit{anti-holomorphic} at places in $\Sigma_{ub}$. Let $F_{h_{ub}}$ denote the function on $\mathfrak{h}^{n}\times\mathbb{C}^{n}$ which is similar to $F_{h}$ but is \textit{skew-holomorphic} (cf. \cite{BS98}, p.80) at places in $\Sigma_{ub}$ (see \S\ref{thethm} for a more precise definition of $F_{h_{ub}}$).

Write the coordinates on $\mathfrak{h}^{n}\times\mathbb{C}^{n}$ as $((\tau_{v})_{v},(z_{v})_{v})$ and define the differential operator
\begin{equation}\label{diffop}
    \Delta_{v}=
    \frac{2i}{\pi}\Big(\frac{\partial}{\partial z_{v}}+4\pi i\frac{z_{v}-\Bar{z}_{v}}{\tau_{v}-\Bar{\tau}_{v}}\Big)
\end{equation}
on the space of smooth functions on $\mathfrak{H}^{n}\times \mathbb{C}^{n}$ for every archimedean place.

Put
\begin{equation*}
    \Delta^{(r)}=\prod_{v\in\Sigma_{\infty}}\Delta_{v}^{r_{v}}.
\end{equation*}

Let $L(s,\Sym^{2}g\times f)$ denote the completed L-function whose central value we are interested in. It is an analytic function on the whole complex plane and satisfies a functional equation of the form
$$L(s,\Sym^{2}g\times f)=\epsilon (s,\Sym^{2}g\times f)L(2\kappa_{0}+2\kappa'_{0}-s,\Sym^{2}g\times f).$$
Here $\kappa_{0}=\max_{v}\{\kappa_{v}\}$ and $\kappa'_{0}=\max_{v}\{\kappa'_{v}\}$. Note that $s=\kappa_{0}+\kappa'_{0}$ is the center of $L(s,\Sym^{2}g\times f)$.

For any cusp forms $f,g$ of weight $\kappa$, let $\langle f,g\rangle$ be the normalized Petersson inner product of $f$ and $g$ defined by
\begin{equation*}
    \langle f,g\rangle=\frac{1}{\Vol(\SL_{2}(\mathcal{O})\backslash \mathfrak{h}^{n})}\displaystyle\int_{\SL_{2}(\mathcal{O})\backslash \mathfrak{h}^{n}}f(x+\imag y)\overline{g(x+\imag y)}\prod_{i=1}^{n}y_{i}^{\kappa_{i}}\frac{dx_{i}dy_{i}}{y_{i}^{2}}
\end{equation*}
The main formula proved in this paper is essentially the following. For a precise version, see Theorem \ref{mainthm}.

\begin{theorem*}
There is a nonzero scalar $c_{\kappa,\kappa'}$, such that
\begin{equation}\label{classver}
    L(\kappa_{0}+\kappa'_{0},\Sym^{2}g\times f)=c_{\kappa,\kappa'}|\langle g_{ub}, \Delta^{(r)}F_{h_{ub}}|_{\mathfrak{h}^{n}}\rangle|^{2}\frac{\langle f,f\rangle}{\langle h_{ub},h_{ub} \rangle},    
\end{equation}
where $\xi_{F}$ is the completed Dedekind zeta function of $F$ and $D_{F}$ is the discriminant of $F$.
\end{theorem*}

\begin{remark}
    The above formula is slightly simplified in the sense that in the actual formula (cf. Theorem \ref{mainthm}) instead of the (non-holomorphic) newform $g_{ub}$ a closely related (non-holomorphic) \textit{oldform} appears. We will also compute $c_{\kappa,\kappa'}$.
\end{remark}

\begin{remark}
Xue (\cite{Xue19}, Prop. 2.1) proved the above formula in the special case of $F=\mathbb{Q}$ and $\kappa,\kappa'$ both odd; a factor of $\xi_{F}(2)^{2}$ appears in our formulation due to the fact that our inner product is normalized.
\end{remark}
We will derive (\ref{classver}) by explicitly computing the local period integrals appearing in the refined Gan-Gross-Prasad formula for $\SL_{2}\times\widetilde{\SL_{2}}$. The main formula will be stated in the language of automorphic representations in \S\ref{mf} and the computations of local integrals will be carried out in \S\ref{archplaces} and \S\ref{nonarchplaces}. As an application of our formula, we will work out the rationality of $L(\kappa_{0}+\kappa'_{0},\Sym^{2}g\times f)$ in the two extreme cases: $\Sigma_{b}=\emptyset$ and $\Sigma_{ub}=\emptyset$. The result on rationality will follow from our main formula by working out the rationality of the factors $\langle g_{ub}, \Delta^{(r)}F_{h_{ub}}|_{\mathfrak{h}^{n}}\rangle$ and $\frac{\langle f,f\rangle}{\langle h_{ub},h_{ub} \rangle}$. This will be done in \S\ref{rationality}. An interesting dichotomy appears when working out the rationality of $\langle g_{ub}, \Delta^{(r)}F_{h_{ub}}|_{\mathfrak{h}^{n}}\rangle$ in the two cases of $\Sigma_{b}=\emptyset$ and $\Sigma_{ub}=\emptyset$: one extreme makes use of Shimura's theory of integral weight (Hilbert) modular forms, the other relies on Shimura's theory of \textit{half-integral} weight modular forms. Our results on rationality of the central L-value are as follows (cf. Theorems \ref{ratb} and \ref{ratub}).

\begin{theorem*}\label{ratclassver}
Put $\epsilon=(\epsilon_{v})\in (\mathbb{Z}/2\mathbb{Z})^{n}\ \text{with}\ \epsilon_{v}=\kappa_{0}+1\ (\text{mod}\ 2)$ at every $v$. Let $u(\cdot,\cdot)$ denote Shimura's $u$-invariant.
\begin{itemize}
    \item [(i)] Suppose $\Sigma_{ub}=\emptyset$. Then for any $\sigma\in \Aut(\mathbb{C})$,
\begin{equation*}
    \sigma\Bigg(\frac{L(\kappa_{0}+\kappa'_{0},\Sym^{2}g\times f)}{\pi^{(\kappa_{0}+2)n-\sum_{i=1}^{n}\kappa_{i}}\sqrt{D_{F}} \langle g,g\rangle^{2}|u(\epsilon,f)|}\Bigg)=\frac{L(\kappa_{0}+\kappa'_{0},\Sym^{2}g^{\sigma}\times f^{\sigma})}{\pi^{(\kappa_{0}+2)n-\sum_{i=1}^{n}\kappa_{i}}\sqrt{D_{F}}\langle g^{\sigma},g^{\sigma}\rangle^{2}|u(\epsilon,f^{\sigma})|},
\end{equation*}

\item [(ii)] Suppose $\Sigma_{b}=\emptyset$. Then for any $\sigma\in \Aut(\mathbb{C}/\tilde{F})$,
\begin{equation*}
    \sigma\Bigg(\frac{L(\kappa_{0}+\kappa'_{0},\Sym^{2}g\times f)}{\pi^{(\kappa_{0}+1)n-\sum_{i=1}^{n}\kappa_{i}}\sqrt{D_{F}} \langle f,f\rangle |u(\epsilon,f)|}\Bigg)=\frac{L(\kappa_{0}+\kappa'_{0},\Sym^{2}g^{\sigma}\times f^{\sigma})}{\pi^{(\kappa_{0}+2)n-\sum_{i=1}^{n}\kappa_{i}}\sqrt{D_{F}}\langle f^{\sigma},f^{\sigma}\rangle |u(\epsilon,f^{\sigma})|}.
\end{equation*}

\end{itemize}
\end{theorem*}

In \S\ref{generalcase}, we will state an explicit conjecture on the rationality of the central L-value in the general case, that is, when both $\Sigma_{b}=\emptyset$ and $\Sigma_{ub}=\emptyset$ can be nonempty. In the rest of the paper we will work in the language of automorphic representations.
\begin{remark}
    It can be verified that the above result on rationality is compatible with Deligne's conjecture on critical values of $L$-functions. In fact, our conjecture (see \S\ref{generalcase}) on rationality in the general case follows by explicitly computing Deligne's period for $L(s,\Sym^{2}g\times f)$.
\end{remark}
\begin{remark}
    In this paper we are assuming full level for both $f$ and $g$. The case of arbitrary levels will be considered in a subsequent paper. In the special case when $F=\mathbb{Q}$ and levels of $f$ and $g$ are square-free, rationality of the central value has been worked out in \cite{CC19}.
\end{remark}

\end{section}

\section{Notations and Preliminaries}\label{Prelim}

\subsection{The base field \texorpdfstring{$F$}{}}$F$ will denote a totally real number field of degree $n$ over $\mathbb{Q}$. Let $\mathcal{O}=\mathcal{O}_{F}$ and $D=D_{F}$ denote, respectively, the ring of integers of $F$ and the discriminant of $F$. Note that $D_{F}$ is a positive integer since $F$ is totally real. We will denote by $\Sigma_{\infty}=\{\alpha_{1},\ldots,\alpha_{n}\}$ the set of archimedean places of $F$, and we put $\alpha=(\alpha_{1},\ldots,\alpha_{n})$ with a fixed order of $\alpha_{j}$. Let $\Sigma_{F}$ denote the set of all places of $F$. Let $\mathbb{A}=\mathbb{A}_{F}$ denote the ring of adeles of $F$. 

For each place $v$ of $F$, let $F_{v}$ denote the completion of $F$ at $v$. Put $F_{\infty}=\prod_{v\in \Sigma_{\infty}}F_{v}$. For each non-archimedean place $v$ of $F$, let $\mathfrak{o}_{v}$ be the ring of integers of $F_{v}$, $\varpi_{v}$ a uniformizer of $\mathfrak{o}_{v}$ and $q_{v}=\# \mathfrak{o}_{v}/ (\varpi_{v})$ the order of the residue field of $F_{v}$.

\subsection{Sum over a subset \texorpdfstring{$S$ of $\Sigma_{\infty}$, $|\cdot|_{S}$}{}}\label{sumsubset}
Let $S\subset \Sigma_{\infty}$. Define
$$\mathbb{Z}^{S}=\{|S|\text{-tuple}\ (a_{v})_{v}\ |\ v\in S\ \text{and}\ a_{v}\in \mathbb{Z}\}.$$

For $\kappa=(\kappa_{v})_{v}\in \mathbb{Z}^{S}$, put
\begin{equation*}
   |\kappa|_{S}=\sum_{v\in S}\kappa_{v}. 
\end{equation*}
Moreover, if $a\in\mathbb{Z}$, then
\begin{equation*}
    |a|_{S}\coloneqq\sum_{v\in S}a=a|S|, 
\end{equation*}
where $|S|$ is the cardinality of the set $S$.
If $a\in\mathbb{Z}$ and $\kappa=(\kappa_{v})_{v}\in \mathbb{Z}^{S}$, then
\begin{equation*}
    |a+\kappa|_{S}\coloneqq |a|_{S}+|\kappa|_{S}.
\end{equation*}
When $S=\Sigma_{\infty}$, we will often denote $|a+\kappa|_{\Sigma_{\infty}}$ simply $|a+\kappa|$.

\subsection{The sets \texorpdfstring{$\Sigma_{ub}$ and $\Sigma_{b}$}{}}
Let $2\kappa=(2\kappa_{v})_{v}$ and $\kappa'+1=(\kappa'_{v}+1)_{v}\in \mathbb{Z}^{\Sigma_{\infty}}$ be weights of some Hilbert modular forms. For any archimedean place $v$, either $(i)\ \kappa_{v}>\kappa'_{v}$ or $(ii)\ \kappa_{v}\leq \kappa'_{v}$. The first case is called the \textit{unbalanced case}, the second case is called the \textit{balanced case}. 
Let $\Sigma_{ub}$ and $\Sigma_{b}$ be subsets of $\Sigma_{\infty}$ whose union is $\Sigma_{\infty}$, such that $\kappa_{v}>\kappa'_{v}$ for $v\in \Sigma_{ub}$ and $\kappa_{v}\leq \kappa'_{v}$ for $v\in \Sigma_{b}$.  Put
\begin{equation*}
    r_{v}=\begin{cases}
    \kappa_{v}-\kappa'_{v}-1 & \text{if}\ v\in \Sigma_{ub}\\
    \kappa'_{v}-\kappa_{v} & \text{if}\ v\in \Sigma_{b}.
    \end{cases}
\end{equation*}

\subsection{The completed Dedekind zeta function \texorpdfstring{$\xi_{F}$}{}}
Let $\xi_{F}$ be the completed Dedekind zeta function of $F$ given by $$\xi_{F}(s)=\Gamma_{\mathbb{R}}(s)^{n}\zeta_{F}(s),$$
where $\Gamma_{\mathbb{R}}(s)=\pi^{-s/2}\Gamma(s/2)$and $\zeta_{F}$ is the Dedekind zeta function of $F$. For Re(s)$>1$, we have $\xi_{F}(s)=\prod_{v}\xi_{F_{v}}(s) $, where
\begin{equation*}
    \xi_{F_{v}}(s)=\Gamma_{\mathbb{R}}(s)\ \ \text{if $v$ is archimedean},
\end{equation*}
and 
\begin{equation*}
    \xi_{F_{v}}(s)=\frac{1}{1-q_{v}^{-s}}\ \ \text{if $v$ is non-archimedean}.
\end{equation*}

\subsection{Measures}\label{nm} We fix the following measures. On $\mathbb{R}$ we take the usual Lebesgue measure. For any finite place $v$ of $F$, let $dx_{v}$ be the Haar measure on $F_{v}$ so that $\Vol(\mathfrak{o}_{v})=1$. Then the measure on $\mathbb{A}$ will be $\prod_{v}dx_{v}$. On $\SL_{2}(\mathbb{R})$, let $dg=y^{-2}dxdyk$, where 
$g=\begin{pmatrix} 
y &  \\
 & y^{-1} 
\end{pmatrix}
\begin{pmatrix} 
1 & x \\
 & 1 
\end{pmatrix}k$ 
is the Iwasawa decomposition and $dk$ is the measure on $\SO_{2}(\mathbb{R})$ so that $\Vol(\SO_{2}(\mathbb{R}))=1$. For any finite place $v$ of $F$, let $dg_{v}$ be the Haar measure on $\SL_{2}(F_{v})$ so that $\Vol(\SL_{2}(\mathfrak{o}_{v}))=1$. Then $D_{F}^{-3/2}\xi_{F}(2)^{-1}\prod_{v}dg_{v}$ is the Tamagawa measure on $\SL_{2}(\mathbb{A})$; that is, the volume of $\SL_{2}(\mathbb{Q})\backslash\SL_{2}(\mathbb{A})$ with respect to this measure is 1. This also gives a measure on $\J(\mathbb{A})$ (cf. \cite{BS98}, Prop. 1.2.4). The inner products will always be defined using these measures, unless otherwise mentioned.

\subsection{The Jacobi group $\J$}
Let $\GL_{4}$ denote the general linear group considered as an algebraic group over $F$. Let $J$ denote the affine algebraic subgroup of $\GL_{4}$ consisting of all matrices which can be written as

\begin{equation*}
   \begin{pmatrix} 
a & & b& \\
 & 1 & & \\
 c& & d &\\
 & & & 1
\end{pmatrix} 
\begin{pmatrix} 
1 & & & \mu\\
 \lambda & 1 & \mu & \xi\\
 & & 1 & -\lambda\\
 & & & 1
\end{pmatrix};\ \ \ \  \begin{pmatrix} 
a & b \\
c & d 
\end{pmatrix}\in \SL_{2}.
\end{equation*}
The algebraic group $\J$ is called the Jacobi group. For brevity, a typical element of the Jacobi group will be written as 
  $ \begin{pmatrix} 
a & b \\
c & d 
\end{pmatrix}
(\lambda,\mu,\xi).$ The subset $\He=\{(\lambda,\mu,\xi)\}\subset\J$ is a subgroup of the Jacobi group called the Heisenberg group. Projection of $\J$ onto $\SL_{2}$ yields a split exact sequence  
$$1\longrightarrow\He\longrightarrow\J\longrightarrow \SL_{2}\longrightarrow 1,$$
so that the Jacobi group becomes the semi-direct product of $\SL_{2}$ and the Heisenberg group. We denote by $B_{\J}$ the Borel subgroup of $\J$ consisting of elements of the form 
$\begin{pmatrix} 
a & b \\
 & d 
\end{pmatrix}
(0,\mu,\xi).$

Let $j=g(X,\kappa)$ and $j'=g'(X',\kappa')$ be two elements in $\J$ with $g,g'\in \SL_{2}$ and $(X,\kappa),(X',\kappa')\in \He $. Then it can be shown that 
\begin{equation*}
    jj'=gg'\Big(Xg'+X',\kappa+\kappa'+\det\begin{pmatrix}
    Xg'\\
    X'
    \end{pmatrix}\Big).
\end{equation*}
A quick calculation shows that the center of the Jacobi group is $Z_{J}=\{(0,0,\kappa)\}$. 

For a local field $k$, the group $\widetilde{\J}(k)=\widetilde{\SL_{2}(k)}\ltimes\He(k)$ will denote the double cover of $\J(k)$, also called the metalplectic Jacobi group. Similarly, for a global field $F$ with ring of adeles $\mathbb{A}$, we have $\widetilde{\J}(\mathbb{A})=\widetilde{\SL_{2}(\mathbb{A})}\ltimes\He(\mathbb{A})$.

\subsection{Lie algebra \texorpdfstring{$\mathfrak{j}_{\mathbb{C}}$}{} of the Jacobi group}\label{liealg}
Let $\mathfrak{j}$ be the Lie algebra of $\J(\mathbb{R})$ and $\mathfrak{j}_{\mathbb{C}}$ its complexification. Then $\mathfrak{j}_{\mathbb{C}}$ is a six dimensional complex vector space spanned by the following elements
\begin{align*}
   \X_{\pm}&=\frac{1}{2}\begin{pmatrix}
1&0 & \pm{i} &0\\
0&0& 0&0\\
\pm{i}&0 &-1 &0\\
0&0&0&0
\end{pmatrix},
&
\Y_{\pm}&=\frac{1}{2}\begin{pmatrix}
0& 0& 0&\pm{i} \\
1&0&\pm{i} &0\\
0&0&0 &-1 \\
0&0&0&0
\end{pmatrix},\\
\Z&=\begin{pmatrix}
0&0 & -i &0\\
0&0& 0&0\\
i& 0&0 &0\\
0&0&0&0
\end{pmatrix},
&
\Z_{0}&=\begin{pmatrix}
0&0 & 0 &0\\
0&0& 0&-i\\
0&0 & 0&0\\
0&0&0&0
\end{pmatrix}.
\end{align*}
Via the identification of $\SL_{2}(\mathbb{R})$ as a subgroup of $\J(\mathbb{R})$, we see that $\mathfrak{sl}_{2,\mathbb{C}}$ is a subalgebra of $\mathfrak{j}_{\mathbb{C}}$ and is spanned by $\{\X_{\pm},\Z\}$. Moreover, $\Z$ spans the (complexified) Lie algebra of the maximal compact subgroup $\SO_{2}(\mathbb{R})$ of $\SL_{2}(\mathbb{R})$. The Lie algebra of the Heisenberg group is spanned by $\{\Y_{\pm},\Z_{0}\}$.

If $\pi$ is a smooth representation of a Lie group and $d\pi$ the induced action of its Lie algebra, then for any Lie algebra element $\X$, $d\pi(\X)v$ will be simply written as $\X v$.

\begin{remark}
At the level of representation theory, $\X_{\pm}$ act as weight $\pm 2$ raising/lowering operators, and $\Y_{\pm}$ as weight $\pm 1$ raising/lowering operators. 
\end{remark}

\subsection{The Schr\"{o}dinger-Weil representation $\omega_{\psi}$}\label{SW rep}
If $k$ is a local field and $\psi$ a non-trivial additive character of $k$, then there is a genuine irreducible representation of $\widetilde{\J(k)}$ known as the Schr\"{o}dinger-Weil representation and denoted by $\omega_{\psi}$. 
It has a realization on the Schwarz space $\mathcal{S}(k)$. Put 
$$m(a)=\begin{pmatrix} 
a &  \\
 & a^{-1} 
\end{pmatrix},\ \ \ \
n(b)=\begin{pmatrix} 
1 & b \\
 & 1 
\end{pmatrix},\ \ \ \
w=\begin{pmatrix} 
 & -1 \\
 1 &  
\end{pmatrix}; \ \ (a\in k^{\times},b\in k).$$
Then the elements $(m(a),1),\ (n(b),1)$, $(w,1)$, $(1,\epsilon)$ and $(\lambda,\mu,\xi)$, for $a,b\in k^{\times}$, $\lambda,\mu,\xi\in k$ and $\epsilon\in \{\pm 1\}$, generate $\widetilde{\J(k)}$. By abuse of notation, we will denote by $g$ an element of the form $(g,1)$ of $\widetilde{\SL_{2}(k)}$.
Then the realization of the Schr\"{o}dinger-Weil representation on the Schwarz space is given by the following formulae.
\begin{equation}\label{weilrep}
\begin{aligned}
  \omega_{\psi}(m(a))\phi(t)&=\frac{\gamma_{\psi}(1)}{\gamma_{\psi}(a)}|a|^{1/2}\phi(at),\\
\omega_{\psi}(n(b))\phi(t)&=\psi(bt^{2})\phi(t),\\
\omega_{\psi}(w)\phi(t)&=\overline{\gamma_{\psi}(1)}|2|^{1/2}\hat{\phi}(-2t),\\
\omega_{\psi}((\lambda,\mu,\xi))\phi(t)&=\psi(\xi+(2t+\lambda)\mu)\phi(x+\lambda)\\
\omega_{\psi}((1,\epsilon))\phi(t)&=\epsilon\phi(t).
\end{aligned}  
\end{equation}

The last formula ensures that the Schr\"{o}dinger-Weil representation is indeed a genuine representation of the metaplectic Jacobi group.
Here $\gamma_{\psi}:k^{\times}\rightarrow \mathbb{C}^{\times}$ denotes the Weil constant and $\hat{\phi}$ denotes the Fourier transform of $\phi$;
\begin{equation*}
       \hat{\phi}(x)=\int_{k}\phi(y)\psi(xy)dy,
  \end{equation*}
where $dy$ is the self-dual additive Haar measure on $k$. The Schr\"{o}dinger-Weil representation is unitary with the inner product 
\begin{equation*}
    \langle\phi_{1},\phi_{2}\rangle=\displaystyle\int_{k}\phi_{1}(y)\overline{\phi_{2}(y)}dy\ \ \ \ \ \ (\phi_{1},\phi_{2}\in \mathcal{S}(k)).
\end{equation*}

For a global field $F$ and a nontrivial additive character $\psi$ of $F\backslash \mathbb{A}$, we have the global Schr\"{o}dinger-Weil representation, $\omega_{\psi}$, of $\widetilde{\J(\mathbb{A})}$,which can be defined by formulae similar to ($\ref{weilrep}$). The representation $\omega_{\psi}$ is unitary with respect to the inner product 
\begin{equation*}
    \langle\phi_{1},\phi_{2}\rangle=\displaystyle\int_{\mathbb{A}}\phi_{1}(y)\overline{\phi_{2}(y)}dy\ \ \ \ \ \ (\phi_{1},\phi_{2}\in \mathcal{S}(\mathbb{A})).
\end{equation*}
Moreover, if $\psi=\otimes_{v}\psi_{v}$, then $\omega_{\psi}=\otimes_{v}\omega_{\psi_{v}}$.

\subsection{The Waldspurger packet $\Wd_{\psi}(\pi)$}\label{wald}

Let $F$ be a number field. Fix a nontrivial additive character $\psi$ of $F$. For each place $v$ of $F$, we consider the pair $(PD^{\times},\widetilde{\SL_{2}})$, where $D$ is a quaternion algebra over $F_{v}$ and $PD^{\times}=D^{\times}/F_{v}^{\times}$. Let $\theta_{PD^{\times}\times\widetilde{\SL_{2}}}(\ \cdot\ ;\psi_{v})$ stand for the theta correspondence with respect to $\psi_{v}$ mapping irreducible representations of $PD^{\times}$ to irreducible representations of $\widetilde{\SL_{2}}(F_{v})$.

Now, let $\pi=\otimes\pi_{v}$ be a cuspidal automorphic representation of $\PGL_{2}(\mathbb{A})$. Put $\sigma_{v}^{+}=\theta_{\PGL_{2}\times\widetilde{\SL_{2}}}(\pi_{v};\psi_{v})$. If $\pi_{v}$ is square integrable, that is, either a discrete series representation if $v$ is archimedean, or, a Steinberg or supercuspidal representation if $v$ is non-archimedean, let $\JL(\pi_{v})$ be the Jacquet-Langlands lift of $\pi_{v}$; it is a representation of $PD^{\times}$, where $D$ is the unique quaternion division algebra over $F_{v}$. Put $\sigma_{v}^{-}=\theta_{PD^{\times}\times\widetilde{\SL_{2}}}(\JL(\pi_{v});\psi_{v})$. The local Waldspurger packet of $\pi_{v}$ with respect to $\psi_{v}$ is given by 
$$\Wd_{\psi_{v}}(\pi_{v})=\begin{cases}
\{\sigma_{v}^{+}\}&\text{if $\pi_{v}$ is not square integrable}\\
\{\sigma_{v}^{+},\sigma_{v}^{-}\} &\text{if $\pi_{v}$ is square integrable}.
\end{cases}$$
Waldspurger (\cite{Wald91}, Prop. 5 and Lemma 40) has shown the following:
\begin{itemize}[noitemsep,topsep=0pt]
    \item [1.]If $F_{v}=\mathbb{R},\ \psi_{v}=\exp(2\pi i \alpha)\ \text{for some real number $\alpha$}$ and $\pi_{v}$ is a discrete series representation of $\PGL_{2}(\mathbb{R})$ of minimal weight $\pm 2\kappa$, then $\sigma_{v}^{+}$ is a discrete series representation of lowest weight $\kappa+\frac{1}{2}$ resp. highest weight $-(\kappa+\frac{1}{2})$ if $\alpha$ is positive resp. negative. Moreover, in this case, $\Wd_{\psi_{v}}(\pi_{v})$ consists of discrete series representations of lowest weight $\kappa+\frac{1}{2}$ and highest weight $-(\kappa+\frac{1}{2})$. 
    \item [2.]For each $\epsilon=(\epsilon_{v})$ such that $\epsilon_{v}$ belongs to $\{\pm1\}$ and takes value 1 when $\pi_{v}$ is not square-integrable, put $\sigma^{\epsilon}=\otimes\sigma_{v}^{\epsilon_{v}}$. Then the global Waldspurger packet of $\pi$ with respect to $\psi$ consists of automorphic representations of $\widetilde{\SL_{2}(\mathbb{A})}$ and is given by 
$$\Wd_{\psi}(\pi)=\Big\{\sigma^{\epsilon} \Big| \prod_{v}\epsilon_{v}=\epsilon\Big(\frac{1}{2},\pi\Big)\Big\}.$$
\end{itemize}

\subsection{The L-function $L(s, \Sym^{2}\tau\times \pi)$}\label{lfct}
Let $\pi$ and $\tau$ be cuspidal automorphic representations of
$\PGL_{2}(\mathbb{A})$ of weights $2\kappa$ and $\kappa'+1$, respectively. For each prime $\mathfrak{p}$, let $\{\alpha_{\mathfrak{p}},\alpha_{\mathfrak{p}}^{-1}\}$ resp. $\{\beta_{\mathfrak{p}},\beta_{\mathfrak{p}}^{-1}\}$ denote the Satake parameters of $\tau$ resp. $\pi$ at $\mathfrak{p}$. 
Put
$$A_{\mathfrak{p}}=\begin{pmatrix}
\alpha_{\mathfrak{p}}^{2} & &\\
&1&\\
&&\alpha_{\mathfrak{p}}^{-2}
\end{pmatrix},\ \ \ \
B_{\mathfrak{p}}=\begin{pmatrix}
\beta_{\mathfrak{p}} &\\
&\beta_{\mathfrak{p}}^{-1}
\end{pmatrix}.$$
Define
$$L_{\fin}(s,\Sym^{2}\tau\times\pi)= \prod_{\mathfrak{p}}\det(I_{6}-(A_{\mathfrak{p}}\otimes B_{\mathfrak{p}})(\N\mathfrak{p})^{-s})^{-1}$$
and (cf. \cite{Jac}, \S 16)
\begin{equation*}\label{linfty}
\begin{split}    
  L_{\infty}(s,\Sym^{2}\tau\times\pi)&=\prod_{v\in \Sigma_{\infty}}\Gamma_{\mathbb{C}}\Big(s+\kappa_{v}-\frac{1}{2}\Big)\Gamma_{\mathbb{C}}\Big(s+\kappa_{v}+\kappa^{'}_{v}-\frac{1}{2}\Big)\\
  &\quad\cdot\prod_{v\in \Sigma_{ub}}\Gamma_{\mathbb{C}}\Big(s+\kappa_{v}-\kappa^{'}_{v}-\frac{1}{2}\Big)
  \cdot\prod_{v\in \Sigma_{b}}\Gamma_{\mathbb{C}}\Big(s+\kappa^{'}_{v}-\kappa_{v}+\frac{1}{2}\Big) 
  \end{split}
\end{equation*}
for Re$(s)>\!>0$. Here $\Gamma_{\mathbb{C}}(s)=2(2\pi)^{-s}\Gamma(s)$. Then the completed L-function $L(s,\Sym^{2}\tau\times\pi)$ is given by 
$$L(s,\Sym^{2}\tau\times\pi)=L_{\infty}(s,\Sym^{2}\tau\times\pi)L_{\fin}(s,\Sym^{2}\tau\times\pi).$$
We say that it is of degree six in $(\N\mathfrak{p})^{-s}$. It has an analytic continuation to the whole complex plane and satisfies the functional equation
$$L(s,\Sym^{2}\tau\times\pi)=\epsilon (s,\Sym^{2}\tau\times \pi)L(1-s,\Sym^{2}\tau\times\pi),$$
where $\epsilon (s,\Sym^{2}\tau\times \pi)$ is a function of exponential type, $ab^{s}$, for some constants $a\in \mathbb{C}^{\times}$ and $b\in \mathbb{R}$, which takes value in $\{\pm 1\}$ when $s=\frac{1}{2}$.

\subsection{Nearly holomorphic automorphic forms}\label{nearhol}
Let $F/\mathbb{Q}$ be a totally real field of degree $n$. Define the Maass-Shimura differential operators $\delta_{v}(c)$ and $\delta_{k}^{(a)}$ on the space of smooth functions on $\mathfrak{h}^{n}$ for $v\in \Sigma_{\infty}$, $c\in \mathbb{R}$, and $0\leq a \in\mathbb{Z}^{n}$ by
\begin{equation*}
    \delta_{v}(c)=(2\pi\imag)^{-1}\Big(\frac{\partial}{\partial z_{v}}+\frac{c}{(z_{v}-\overline{z_{v}})}\Big),
\end{equation*}
\begin{equation*}
    \delta_{\kappa}^{(a)}=\displaystyle\prod_{v\in\Sigma_{\infty}}\displaystyle\prod_{j=1}^{a_{v}}\delta_{v}(\kappa_{v}+2j-2).
\end{equation*}
By a nearly holomorphic Hilbert modular form of (integral or half-integral) weight $\kappa$ and level $\mathfrak{n}$, we mean (cf. \cite{Shi87}, Lemma 8.2) a real analytic function $f$ on $\mathfrak{h}^{n}$ such that
\begin{equation*}
    f=\sum_{0\leq p\leq\kappa/2}\delta_{\kappa-2p}^{(p)}g_{p}+\begin{cases}
    c\delta_{2}^{(\kappa/2)-1}E_{2}\ \ \ \text{if}\ F=\mathbb{Q}\  \text{and}\ \kappa\in 2\mathbb{Z}\\
    0\ \ \ \ \ \ \text{otherwise}\end{cases},
\end{equation*}
where $g_{p}$ is a Hilbert modular form of weight $\kappa-2p$ and level $\mathfrak{n}$, $c\in \mathbb{C}$, and
$$E_{2}(z)=(4\pi y)^{-1}-12^{-1}+2\displaystyle\sum_{n=1}^{\infty}\Big(\displaystyle\sum_{0<d|n}d\Big)e^{2\pi \imag nz}.$$ 
Note that, $\delta_{\kappa}\coloneqq\delta_{\kappa}^{(1)}$ acts as a weight raising operator -- raising the weight by 2 --
on the space of nearly holomorphic modular forms.

Next, we will define a nearly holomorphic modular form in the language of automorphic forms. Since a function on $\mathfrak{h}^{n}$ can be realized as a function on $G=\GL_{2}^{+}(\mathbb{R})^{n}$, we may think of $\delta_{\kappa}^{(a)}$ as a differential operator on the space of smooth functions on $G$. It is not hard to check that $\delta_{\kappa}^{(a)}$ is in fact a left-invariant differential operator and hence, under the identification of left-invariant differential operators on $G$ with he universal enveloping algebra $\mathcal{U}(\Lie(G))$, may be viewed as an element $\mathbf{\delta}_{\kappa}^{(a)}\in \mathcal{U}(\Lie(G))$. Recall also the element $\X_{+}$ of $\Lie(\J(\mathbb{R}))$ (see \S\ref{liealg}) which as an element of $\Lie(\GL_{2}^{+}(\mathbb{R}))$ looks like
$$\X_{+}=\frac{1}{2}\begin{pmatrix}
1 & \imag\\
\imag &-1
\end{pmatrix}.$$
It is a fact that $\X_{+}$ as a differential operator on $\GL_{2}^{+}(\mathbb{R})$ is a weight raising operator, raising the weight by 2. It follows from the structure of irreducible $(\Lie(\GL_{2}^{+}(\mathbb{R})),SO_{2}(\mathbb{R}))$-modules (cf. p. 24) that if $g_{p}$ is a Hilbert modular form of weight $\kappa-2p$ and $\mathbf{g}_{p}$ denotes its adelization, then $\mathbf{\delta}_{\kappa-2p}^{(p)}\mathbf{g}_{p}$ must be a scalar multiple of $\displaystyle\prod_{v\in\Sigma_{\infty}}\X_{+}^{p_{v}}\mathbf{g}_{p}$. Hence, we may call an automorphic form $\varphi$ a nearly holomorphic automorphic form of weight $\kappa$ if (assume $F\neq \mathbb{Q}$)
\begin{equation*}
    \varphi=\displaystyle\sum_{0\leq p\leq\kappa/2}\displaystyle\prod_{v\in\Sigma_{\infty}}\X_{+}^{p_{v}}\mathbf{g}_{p},
\end{equation*}
where $\mathbf{g}_{p}$ is a Hilbert modular automorphic form of weight $\kappa-2p$. 

\begin{lemma}\label{nearhollemma}
Let $\varphi$ be an automorphic form on $\GL_{2}(\mathbb{A})$. Then $\varphi$ is a linear combination of nearly holomorphic automorphic forms if and only if $\displaystyle\prod_{v\in\Sigma_{\infty}}\X_{-}^{s_{v}}\varphi=0$ for some $s=(s_{v})_{v}\in\mathbb{Z}^{\Sigma_{\infty}}$, where $\X_{-}$ is the element of $\mathfrak{j}_{\mathbb{C}}$ defined in $\S\ref{liealg}$ but viewed as an element of $\Lie(\SL_{2})$, the (complexified) Lie algebra of $\SL_{2}$.
\end{lemma}
\begin{proof}
Without loss of generality assume that $\varphi$ is nearly holomorphic of the form $\varphi=\prod_{v\in\Sigma_{\infty}}\X_{+}^{p_{v}}\mathbf{g}_{p}$, where $\mathbf{g}_{p}$ is a Hilbert modular automorphic form of weight $\kappa-2p$. It follows from the theory of $(\Lie(\GL_{2}(\mathbb{R})),O_{2}(\mathbb{R}))$-modules that $\X_{-}$ acts as a weight lowering operator, lowering the weight by 2, and that $\X_{-}\X_{+}$ acts as a scalar (cf. \cite{Bump}, Prop. 2.5.2). Therefore, $\prod_{v\in\Sigma_{\infty}}\X_{-}^{p_{v}+1}\varphi=0$.

For the converse, write $\varphi$ as a linear combination of linearly independent pure tensors $\phi_{i}$. If $\prod_{v\in\Sigma_{\infty}}\X_{-}^{s_{v}}\varphi=0$, then linear independence implies that $\prod_{v\in\Sigma_{\infty}}\X_{-}^{s_{v}}\varphi_{i}=0$ for all $i$. Hence, without loss of generality, we may assume $\varphi$ is a pure tensor. If $\prod_{v\in\Sigma_{\infty}}\X_{-}^{s_{v}}\varphi=0$ and $s=(s_{v})_{v}$ is smallest such, then $\mathbf{f}=\prod_{v\in\Sigma_{\infty}}\X_{-}^{s_{v}-1}\varphi$ must be a holomorphic Hilbert modular automorphic form. It then follows from the theory of $(\Lie(\GL_{2}(\mathbb{R})),O_{2}(\mathbb{R}))$-modules that
$\prod_{v\in\Sigma_{\infty}}\X_{+}^{s_{v}-1}\mathbf{f}=c\phi\ \ \ \text{for some\ $c\in \mathbb{C}$.}$
\end{proof}

\begin{section}{The Main Formula}\label{mf}
Let $\pi=\otimes\pi_{v}$ resp. $\tau=\otimes\tau_{v}$ be an irreducible cuspidal automorphic representations of $\PGL_{2}(\mathbb{A})$ of weights $2\kappa$ resp. $\kappa'+1$. Let $\psi$ be the additive character of $F\backslash \mathbb{A}$ whose infinity component is given by $x\mapsto \exp(2\pi ix)$ for any real place. For any non-archimedean place of $F$, let $\delta_{v}$ denote the conductor of $\psi_{v}$, that is, $\delta_{v}^{-1}\mathfrak{o}_{v}$ is the largest subgroup of $F_{v}$ on which $\psi_{v}$ is trivial. We first make an observation.
\begin{lemma} $L(\frac{1}{2},\Sym^{2}\tau\times \pi)=0$, unless $|\kappa|\equiv |\Sigma_{b}|\ \ (\text{mod}\ 2)$.
\end{lemma}
\begin{proof}
By (\cite{HK91}, \S8),
\begin{equation*}
   \epsilon_{v}\Big(\frac{1}{2},\tau\times \tau\times \pi\Big)=\begin{cases}
1 & \text{if}\ v\in \Sigma_{ub}\\
-1 & \text{if}\ v\in \Sigma_{b}.
\end{cases} 
\end{equation*}
Also, since $\pi$ is of full level, $\epsilon(\frac{1}{2},\pi)=(-1)^{|\kappa|}$. 
Since
\begin{equation*}
   \epsilon\Big(\frac{1}{2},\tau\times \tau\times \pi\Big)=\epsilon\Big(\frac{1}{2},\Sym^{2}\tau\times \pi\Big)\epsilon\Big(\frac{1}{2},\pi\Big),
\end{equation*}
we see that if $|\kappa|$ is odd (resp. even) and $|\Sigma_{b}|$ is even (resp. odd), then $\epsilon (\frac{1}{2},\Sym^{2}\tau\times \pi)=-1$ and thus $L(\frac{1}{2},\Sym^{2}\tau\times \pi)=0$.
\end{proof} Hence, from now on we will assume that
$$|\kappa|\equiv |\Sigma_{b}|\ \ (\text{mod}\ 2). \ \ (\star) $$

Next, we will choose an automorphic representation $\sigma$ in the representations $\Wd_{\overline{\psi}}(\pi)$ and vectors in $\tau,\ \sigma,\ \omega_{\psi}$ so that the local integrals appearing in the refined Gan-Gross-Prasad conjecture for $\SL_{2}\times\widetilde{\SL_{2}}$ are nonzero. We will define these by specifying their local components.

\subsection{Local Components}\label{loccomp}

Suppose that $v\in \Sigma_{b}$.
\begin{itemize}[noitemsep,topsep=0pt]
    \item $\pi_{v}$ is a discrete series representation of $\PGL_{2}(\mathbb{R})$ of minimal weight $\pm 2\kappa_{v}$. Let $\sigma_{v}$ be the discrete series representation of $\widetilde{\SL_{2}(\mathbb{R})} $ of lowest weight $\kappa_{v}+\frac{1}{2}$ and $h_{v}$ denote a lowest weight vector of $\sigma_{v}$.
    \item $\tau_{v}$ is a discrete series representation of $\PGL_{2}(\mathbb{R})$ of minimal weight $\pm(\kappa'_{v}+1)$. Let $g_{v}\in \tau_{v}$ be a lowest weight vector.   
\end{itemize}
Suppose that $v\in \Sigma_{ub}$.
\begin{itemize}[noitemsep,topsep=0pt]
    \item $\pi_{v}$ is a discrete series representation of $\PGL_{2}(\mathbb{R})$ of minimal weight $\pm 2\kappa_{v}$. Consider the discrete series representation of $\widetilde{\SL_{2}(\mathbb{R})}$ of lowest weight $(\kappa_{v}+\frac{1}{2})$ and let $h_{v}$ denote a lowest weight vector. We take $\sigma_{v}$ to be the discrete series representation of \textit{highest weight} $-(\kappa_{v}+\frac{1}{2})$. The complex conjugate $\overline{h_{v}}$ is a highest weight vector of $\sigma_{v}$.
    \item $\tau_{v}$ is a discrete series representation of $\PGL_{2}(\mathbb{R})$ of mimimal weight $\pm(\kappa'_{v}+1)$. Let $g_{v}\in \tau_{v}$ be a lowest weight vector. Then its complex conjugate, $\overline{g_{v}}$, is a highest weight vector in $\tau_{v}$.
    \end{itemize}
Suppose that $v<\infty$.
\begin{itemize}[noitemsep,topsep=0pt]
    \item $\pi_{v}$ is an unramified principal series representation of $\PGL_{2}(F_{v})$. The representation $\sigma_{v}$ will be defined precisely in \S7. For now, we mention that it is an unramified principal series representation of $\widetilde{\SL_{2}(F_{v})}$ of the form $\widetilde{I}_{\overline{\psi}}(s)$, and contains a distinguished vector, $h_{v}$, fixed by the idempotent $E^{K}$ of the Hecke algebra (with respect to $\overline{\psi}$), also defined in \S7. Put $h_{v}^{(2)}=\sigma(m(2))h_{v}$.
    \item $\tau_{v}$ is an unramified principal series representation of $\PGL_{2}(F_{v})$. Let $g_{v}\in \tau_{v}$ be a  $\PGL_{2}(\mathcal{O}_{v})$-fixed vector. Let
$g_{v}^{(\delta_{v})}=\tau\begin{pmatrix} 
\delta_{v}^{-1} &  \\
 & 1 
\end{pmatrix} g_{v}.$
\end{itemize}

Put
\begin{equation*}
    \sigma=\otimes\sigma_{v},\ \ \ \ \mathbf{g}^{(\delta)}_{ub}=\otimes_{v\in\Sigma_{ub}}\overline{g_{v}}\otimes_{v\in\Sigma_{b}}g_{v}\otimes_{v<\infty}g_{v}^{(\delta_{v})},\ \ \ \ \ \mathbf{h}_{ub}^{(2)}=\otimes_{v\in\Sigma_{ub}}\overline{h_{v}}\otimes_{v\in\Sigma_{b}}h_{v}\otimes_{v<\infty}h_{v}^{(2)}.
\end{equation*}
Note that, by our assumption $(\star)$, $\sigma$ indeed lies in  $\Wd_{\overline{\psi}}(\pi)$.  In particular, it is an automorphic representation of $\widetilde{\SL_{2}(\mathbb{A}})$. Put $\mathbf{h}_{ub}=\otimes_{v\in\Sigma_{ub}}\overline{h_{v}}\otimes_{v\notin\Sigma_{ub}}h_{v}$. Clearly, $\sigma$ is generated by $\mathbf{h}_{ub}$. Moreover, $\mathbf{h}^{(2)}_{ub}=\sigma(m(2_{\fin}))\mathbf{h}_{ub}$, where $2_{\fin}$ denotes the finite part of the principal idele $2\in F^{\times}$. 
\begin{remark}
$\mathbf{h}_{ub}$ is the adelization of the function $h_{ub}$ defined in the Introduction. Moreover, the vector $\mathbf{h}^{(2)}_{ub}$ is (up to a nonzero constant number) the adelization of the function $h_{ub}(\tau/4)$, where $\tau$ is the coordinate on $\mathfrak{h}^{n}$.
\end{remark}

We now fix a vector in the Schr\"{o}dinger-Weil representation $\omega_{\psi}=\otimes \omega_{\psi_{v}}$. Let $\boldsymbol{\phi}=\otimes\phi_{v}\in \mathcal{S}(\mathbb{A})$ be the Schwarz function with 
$$\phi_{v}=\begin{cases}
\exp(-2\pi  x^{2}) & \text{if}\ v\in \Sigma_{\infty}\\
\mathds{1}_{\frac{1}{2}\mathfrak{o}_{v}} & \text{if}\ v|2\\
\mathds{1}_{\mathfrak{o}_{v}} & \text{otherwise}
\end{cases}$$
Recall that \begin{equation*}
    r_{v}=\begin{cases}
    \kappa_{v}-\kappa'_{v}-1 & \text{if}\ v\in \Sigma_{ub}\\
    \kappa'_{v}-\kappa_{v} & \text{if}\ v\in \Sigma_{b}.
    \end{cases}
\end{equation*}
Set 
\begin{equation*}
    Y_{+}^{(r)}\boldsymbol{\phi}=(\otimes_{v\in \Sigma_{\infty}}Y_{+}^{r_{v}}\phi_{v})\otimes(\otimes_{v<\infty}\phi_{v}),
\end{equation*}
where $Y_{+}=\Big(\frac{2}{\pi}\Big)\cdot\Y_{+}$ is the normalization of the element $\Y_{+}$ of $\mathfrak{j}_{\mathbb{C}}$, the complexified Lie algebra of the Jacobi group, defined in \S\ref{liealg}, and acts as a weight raising differential operator, raising the weight by 1.

\begin{remark}
  $Y_{+}^{(r)}\in \mathfrak{j}_{\mathbb{C}}$ is the adelization of the differential operator $\Delta^{(r)}$ defined in \S\ref{classlang}.  
\end{remark}

\begin{lemma}\label{normphi}
We have,
\begin{equation}
    \langle Y_{+}^{(r)}\boldsymbol{\phi},Y_{+}^{(r)}\boldsymbol{\phi}\rangle=\prod_{v\in \Sigma_{\infty}}\frac{(2r_{v})!}{r_{v}!\pi^{r_{v}}}.
\end{equation}
\end{lemma}
\begin{proof}
For a non-archimedean place $v$ of odd residue characteristic it is not hard to see that $\langle\phi_{v},\phi_{v}\rangle=1$. Moreover, if $v|2$, then $\langle\phi_{v},\phi_{v}\rangle=2^{d}$, since $[\frac{1}{2}\mathfrak{o}_{v}:\mathfrak{o}_{v}]=2^{d}$ and $\Vol(\mathfrak{o}_{v})=1$. Let $v\in \Sigma_{\infty}$. Then by Lemma 3.2.1 in \cite{BS98}, we have
\begin{align*}
    Y_{+}^{r_{v}}\exp(-2\pi  x^{2})&=\Big(\frac{2}{\pi}\Big)^{r_{v}}\Big(\frac{1}{2}\frac{d}{dx}-2\pi  x\Big)^{r_{v}}\exp(-2\pi  x^{2})\\
    &=\Big(\frac{2}{\pi}\Big)^{r_{v}}\cdot (-1)^{r_{v}}\Big(\frac{\pi }{2}\Big)^{r_{v}/2}H_{r_{v}}(\sqrt{2\pi }x)\exp(-2\pi  x^{2}),
\end{align*}
where $H_{r}$ is the $r$th Hermite polynomial. It is a fact that (cf. \cite{GR00}, \S7.375, \#1)
\begin{equation*}
    \displaystyle\int_{-\infty}^{\infty}H_{r}(y)H_{r}(y)e^{-2y^{2}}dy=\frac{(2r)!}{2^{r}r!}\sqrt{\frac{\pi}{2}}.
\end{equation*}
Therefore 
\begin{equation*}
    \langle Y_{+}^{r_{v}}\phi_{v},Y_{+}^{r_{v}}\phi_{v}\rangle=\frac{(2r_{v})!}{2\pi^{r_{v}}r_{v}!}.
\end{equation*}
Since $\langle Y_{+}^{(r)}\boldsymbol{\phi},Y_{+}^{(r)}\boldsymbol{\phi}\rangle=\prod_{v|\infty}\langle Y_{+}^{r_{v}}\phi_{v},Y_{+}^{r_{v}}\phi_{v}\rangle\prod_{v<\infty}\langle\phi_{v},\phi_{v}\rangle$, the lemma follows .
\end{proof}

\subsection{The Theorem}\label{thethm}
To state concisely the next couple of propositions, we introduce some new notation. Let 
\begin{equation*}
    (g_{v}',h_{v}',\phi_{v}')=
    \begin{cases}
    (g_{v},h_{v},Y_{+}^{r_{v}}\phi_{v}) \ \ \ \ \text{if}\ v\in \Sigma_{b}\\
    (\overline{g_{v}},\overline{h_{v}},Y_{+}^{r_{v}}\phi_{v}) \ \ \ \ \text{if}\ v\in \Sigma_{ub}\\
    (g_{v}^{(\delta_{v})},h_{v}^{(2)},\phi_{v}) \ \ \ \ \ \ \text{if}\ v<\infty
    
    \end{cases}
\end{equation*}
For each place $v$ pick some inner product on $\sigma_{v}, \tau_{v}$ and $\omega_{\psi_{v}}$. Define 
$$ \Pe_{v}( g_{v}',h_{v}',\phi_{v}')=\int_{\SL_{2}(F_{v})}\langle \tau_{v}(\alpha_{v})g_{v}',g_{v}'\rangle\overline{\langle \sigma_{v}(\alpha_{v})h_{v}',h_{v}'\rangle \langle \omega_{\psi_{v}}(\alpha_{v})\phi_{v}',\phi_{v}'\rangle}  d\alpha_{v}.
$$
By Lemma 4.3 in \cite{Qiu14}, $\Pe_{v}$ is a finite number. We normalize it by setting
$$\Pe_{v}^{\sharp}(g_{v}',h_{v}',\phi_{v}')=\Bigg(\frac{\xi_{F_{v}}(2)L(1/2,\Sym^{2}\tau_{v}\times \pi_{v})}{L(1,\pi_{v},\Ad)L(1,\tau_{v},\Ad)}\Bigg)^{-1}\frac{\Pe_{v}(g_{v}',h_{v}',\phi_{v}')}{ \langle g_{v}',g_{v}'\rangle \langle h_{v}',h_{v}'\rangle \langle \phi_{v}',\phi_{v}'\rangle}.$$
Note that in the notation introduced above, 
\begin{equation*}
    \mathbf{g}^{(\delta)}_{ub}=\otimes_{v}g_{v}',\ \ \ \mathbf{h}_{ub}^{(2)}=\otimes_{v}h_{v}',\ \ \ Y_{+}^{(r)}\boldsymbol{\phi}=\otimes_{v}\phi_{v}'.
\end{equation*}

We will also need the automorphic version of the Schr\"{o}dinger-Weil representation $\omega_{\psi}$. Let $\boldsymbol{\phi}$ be any Schwarz function. To it one can associate a function on $\widetilde{\J(\mathbb{A})}$ defined as $\Theta_{\boldsymbol{\phi}}(g)=\sum_{a\in F}(\omega_{\psi}(g))Y_{+}^{(r)}\boldsymbol{\phi}(a)$ which is known to be an automorphic form on $\widetilde{\J(\mathbb{A})}$. The action of $\widetilde{\J(\mathbb{A})}$ on the space of $\Theta_{\boldsymbol{\phi}}$ by right translation then defines an automorphic representation of $\widetilde{\J(\mathbb{A})}$ which is isomorphic to the representation $\omega_{\psi}$ (cf \cite{BS98}, \S7.2). 

Now, let $F_{\mathbf{h}_{ub}}$ denote the automorphic form on $\J(\mathbb{A})$ defined by $F_{\mathbf{h}_{ub}}(gh)=\mathbf{h}^{(2)}_{ub}(g)\Theta_{\boldsymbol{\phi}}(gh),$ where $g\in SL_{2}(\mathbb{A})$, $h\in \He(\mathbb{A})$. The automorphic representation generated by $F_{\mathbf{h}_{ub}}$ is isomorphic to $\sigma\otimes \omega_{\psi}$ (by \cite{BS98}, 7.3.3). Therefore, $Y_{+}^{(r)}F_{\mathbf{h}_{ub}}(gh)=\mathbf{h}^{(2)}_{ub}(g)\Theta_{Y_{+}^{(r)}\boldsymbol{\phi}}(gh)$.

For any $f\in L^{2}(\SL_{2}(F)\backslash\SL_{2}(\mathbb{A}))$, let $\langle f,f\rangle$ be the Petersson norm of $f$ defined by 
\begin{equation*}
    \langle f,f\rangle=\displaystyle\int_{\SL_{2}(F)\backslash\SL_{2}(\mathbb{A})}|f(X)|^{2}dX,
\end{equation*}
where $dX$ is the Tamagawa measure on $\SL_{2}(\mathbb{A})$.
\begin{prop}\label{refinedGGP}
For $\mathbf{g}^{(\delta)}_{ub}\in\tau,\mathbf{h}_{ub}^{(2)}\in\sigma$ and $Y_{+}^{(r)}\boldsymbol{\phi}\in\omega_{\psi}$ as above,
\begin{equation*}
  \resizebox{0.99\hsize}{!}{$\frac{\big|\langle \mathbf{g}^{(\delta)}_{ub}, (Y_{+}^{(r)}F_{\mathbf{h}_{ub}})\big|_{\SL_{2}(\mathbb{A})}\rangle\big|^{2} }{\langle\mathbf{g}^{(\delta)}_{ub},\mathbf{g}^{(\delta)}_{ub}\rangle\langle\mathbf{h}_{ub}^{(2)},\mathbf{h}_{ub}^{(2)}\rangle\langle Y_{+}^{(r)}\boldsymbol{\phi},Y_{+}^{(r)}\boldsymbol{\phi}\rangle}=\frac{D_{F}^{-3/2}}{4}\times \frac{L(\frac{1}{2},\Sym^{2}\tau\times \pi)}{L(1,\pi,\Ad)L(1,\tau,\Ad)}\times
  \prod_{v\in \Sigma_{F}}\Pe_{v}^{\sharp}( g_{v}', h_{v}',\mathbf{\phi}_{v}')$},
\end{equation*}
where $\langle-,-\rangle$ stands for the standard inner product on $L^{2}(\SL_{2}(F)\backslash\SL_{2}(\mathbb{A}))$.
\end{prop}

This proposition is Theorem 4.5 in \cite{Qiu14} wherein we have taken the factorizable vectors to be $\mathbf{g}^{(\delta)}_{ub},\ \mathbf{h}_{ub}^{(2)}$ and $Y_{+}^{(r)}\boldsymbol{\phi}$. The above is also referred to as the refined Gan-Gross-Prasad formula for $\SL_{2}\times\widetilde{\SL_{2}}$. 

\begin{remark}
    In the above proposition, the automorphic form $\mathbf{g}^{(\delta)}_{ub}$ on $\PGL_{2}(\mathbb{A})$ is being viewed as an automorphic form on $\SL_{2}(\mathbb{A})$ via the surjection $\SL_{2}(\mathbb{A})\rightarrow \GL_{2}(\mathbb{A})\rightarrow\PGL_{2}(\mathbb{A})$, where the second map is the quotient map and the first is the inclusion map.
\end{remark}

\begin{remark}
Comparing the formula in the above proposition with that in Theorem 4.5 of \cite{Qiu14}, the first factor on the right is different due to the different choice of measures. 
\end{remark}

In \S\ref{archplaces} and \S\ref{nonarchplaces} we'll show that,
\begin{prop}\label{localperiod}
 $$\Pe_{v}^{\sharp}( g_{v}', h_{v}',\mathbf{\phi}_{v}')=
 \begin{cases}
 \frac{2^{2r_{v}+1}}{r_{v}!}\pi^{r_{v}} & \text{if}\ v\in \Sigma_{ub}\\
 \frac{1}{r_{v}!}\pi^{r_{v}} & \text{if}\ v\in \Sigma_{b}\\
 2^{-[F_{v}:\mathbb{Q}_{2}]} & \text{if}\ v|2\\
 1 &\text{otherwise}
 \end{cases}
 $$
 
\end{prop}
We can now state and prove our main formula which is an explicit version of Proposition \ref{refinedGGP}. Note that, even though Proposition \ref{refinedGGP} is independent of any particular normalization of  $\mathbf{g}^{(\delta)}_{ub}$ (or for that matter, $\mathbf{h}^{(2)}_{ub}$ or $Y_{+}^{(r)}F_{\mathbf{h}_{ub}}$), to state the explicit version we normalize it as follows. Put $\mathbf{g}=\otimes_{v\in\Sigma_{F}}g_{v}$, where $g_{v}$'s are as defined in \S\ref{loccomp}. Clearly, $\mathbf{g}$ is a holomorphic newform. Normalize $\mathbf{g}$ so that,
\begin{equation*}
    W_{\mathbf{g}}\begin{pmatrix}
        \eta^{-1} & \\
        & 1
    \end{pmatrix}
    =e^{-2\pi n},
\end{equation*}
where $W_{\mathbf{g}}$ is the Whittaker function of $\mathbf{g}$ defined by 
\begin{equation*}
    W_{\mathbf{g}}(X)=\displaystyle\int_{F\backslash\mathbb{A}}\mathbf{g}\bigg(\begin{pmatrix}
        1& y\\
        & 1
    \end{pmatrix}X\bigg)\overline{\psi(y)}dy
\end{equation*}
with the Tamagawa measure $dy$ on $\mathbb{A}$, $\eta=(\varpi_{v}^{n_{v}})\in \mathbb{A}^{\times}_{\fin}$, and $n=[F:\mathbb{Q}]$. We use the same normalization for $\mathbf{g}^{(\delta)}_{ub}$. 
\begin{theorem}\label{mainthm}
We have,
\begin{equation}\label{mainformula}
    L\Big(\frac{1}{2},\Sym^{2}\tau\times \pi\Big)=\frac{2^{a}D_{F}^{2}}{\displaystyle\prod_{v\in \Sigma_{\infty}}\binom{2r_{v}}{r_{v}}}\xi_{F}(2)|\langle \mathbf{g}^{(\delta)}_{ub}, (Y_{+}^{(r)}F_{\mathbf{h}_{ub}})|_{\SL_{2}(\mathbb{A})}\rangle|^{2}\frac{L(1,\pi,\Ad)}{\langle\mathbf{h}_{ub},\mathbf{h}_{ub}\rangle},
\end{equation}
where $a=|\kappa'+4|_{\Sigma_{\infty}}+|2\kappa'-2\kappa+1|_{\Sigma_{ub}}+2$.
\end{theorem}
\begin{proof}
It is a fact that (cf. \cite{IP16}, Prop. 6.6, Lemma 6.1)
\begin{align}\label{adjLvalue}
    L(1,\tau,\Ad)=2^{2n+|\kappa'+1|}D_{F}^{1/2}\xi_{F}(2)\langle\mathbf{g},\mathbf{g}\rangle.
\end{align}
Moreover, by comparing the local components of $\mathbf{g}$ and $\mathbf{g}^{(\delta)}_{ub}$, it is not hard to see that the local inner products, $\langle\mathbf{g}^{(\delta)}_{ub},\mathbf{g}^{(\delta)}_{ub}\rangle_{v}$ and $\langle\mathbf{g},\mathbf{g}\rangle_{v}$, are equal at all places; thus $\langle\mathbf{g}^{(\delta)}_{ub},\mathbf{g}^{(\delta)}_{ub}\rangle=\langle\mathbf{g},\mathbf{g}\rangle$. Furthermore, since $\mathbf{h}^{(2)}_{ub}=\sigma(m(2_{\fin}))\mathbf{h}_{ub}$ and $\sigma$ is a unitary representation, therefore $\langle\mathbf{h}^{(2)}_{ub},\mathbf{h}^{(2)}_{ub}\rangle=\langle\mathbf{h}_{ub},\mathbf{h}_{ub}\rangle$. The theorem now follows by combining Proposition \ref{localperiod} and Lemma \ref{normphi} with Proposition \ref{refinedGGP}.
\end{proof}

\begin{cor}
$L\Big(\frac{1}{2},\Sym^{2}\tau\times \pi\Big)$ is a nonnegative real number.
\end{cor}
\end{section}

\begin{remark}
    Theorem \ref{mainthm} is the precise formulation in adelic language of the formula (\ref{classver}) stated in the Introduction. We mention that, in the true formula, the (non-holomorphic) newform $g_{ub}$ gets replaced by the closely related (non-holomorphic) \textit{oldform} $\mathbf{g}^{(\delta)}_{ub}$, and, instead of $\langle f,f \rangle$, we now have $L(1,\pi, \Ad)$, both values related to each other by a formula similar to (\ref{adjLvalue}).
\end{remark}

\begin{section}{Rationality Of The Central Value In Some Special Cases}\label{rationality}
To work out the rationality of $L(\frac{1}{2},\Sym^{2}\tau\times\pi)$, we need to work out the rationality of the factors $|\langle \mathbf{g}^{(\delta)}_{ub}, (Y_{+}^{(r)}F_{\mathbf{h}_{ub}})|_{\SL_{2}(\mathbb{A})}\rangle|^{2}$ and $\frac{L(1,\pi,\Ad)}{\langle\mathbf{h}_{ub},\mathbf{h}_{ub}\rangle}$ appearing in our formula (\ref{mainformula}). 

Rationality of the second factor follows from the generalization of the Kohnen-Zagier formula proved in \cite{HI13}. We state this formula here following \cite{HI13}. Let $\eta\in\mathcal{O}^{\times}$ such that $\N_{F/\mathbb{Q}}\eta=(-1)^{|\kappa|}$. Then there exists a totally positive $\xi\in F^{\times}$ such that $L(\frac{1}{2},\pi \otimes\chi_{\eta \xi})\neq 0$, where $L(s,\pi \otimes\chi_{\eta \xi})$ is the completed L-function of $\pi \otimes\chi_{\eta \xi}$. Put $\sigma'=\theta(\pi\otimes \chi_{\eta\xi},\psi_{\xi})$, where $\psi_{\xi}(x)\coloneqq \psi(\xi x)$, and let $\mathbf{h}\in\sigma'$ such that $\mathbf{h}$ lies in the Kohnen plus space, $S_{\kappa+\frac{1}{2}}(\Gamma)^{+}=S_{\kappa+\frac{1}{2}}(\Gamma)^{E^{K}}$ (cf. \cite{HI13} for notation). In particular, $\mathbf{h}$ is a holomorphic vector of weight $\kappa+\frac{1}{2}$ and is fixed by the idempotent $E^{K}$ of the Hecke algebra of $\widetilde{\SL_{2}}$ with respect to the character $\psi_{\eta}$. By scaling $\mathbf{h}$, we may assume that $\langle\mathbf{h},\mathbf{h}\rangle=\langle\mathbf{h}_{ub},\mathbf{h}_{ub}\rangle$, where $\mathbf{h}_{ub}$ is as defined in \S\ref{loccomp}.

Let
$\sum_{\xi'\in\mathcal{O}}c(\xi')e^{2\pi i\xi' z}$ be the Fourier expansion of $\mathbf{h}$ (or more precisely, of the half-integral weight modular form  whose adelization is $\mathbf{h}$). Then we can assume that $c(\xi)\neq 0$ (we may have to replace $\xi$ by $\xi a^{2}$, for some $a\in F^{\times}$, to ensure this; see Lemma 12.4 in \cite{HI13}). Consider the following normalization of $\mathbf{h}$. Put $\mathbf{h}\coloneqq \mathbf{h}\frac{|\Psi(\eta\xi,\pi)|}{|c(\xi)|}$, where $\Psi(\eta\xi,\pi)\mspace{-4mu}=\mspace{-4mu}\prod_{v<\infty}\Psi_{v}(\eta\xi,\alpha_{v})$ is as defined in \cite{HI13}. We apply the same normalization to $\mathbf{h}_{ub}$; thus, $\langle\mathbf{h},\mathbf{h}\rangle=\langle\mathbf{h}_{ub},\mathbf{h}_{ub}\rangle$ still holds. Then, by Eqn. (12.1) in \cite{HI13}, it follows that, 

\begin{lemma}\label{cor3.3}
\begin{equation*}
    \frac{L(1,\pi,\Ad)}{\langle\mathbf{h}_{ub},\mathbf{h}_{ub}\rangle}=D_{F}^{1/2}2^{-1+3|\kappa|}\xi_{F}(2)L(\frac{1}{2},\pi \otimes\chi_{\eta \xi})\prod_{v\in\Sigma_{\infty}}\xi_{v}^{\kappa_{v}-\frac{1}{2}},
\end{equation*}
where $\chi_{\eta\xi}$ is the character attached to $F(\sqrt{\eta\xi})/F$. 
\end{lemma}

\begin{remark}\label{rem_nor}
Henceforth, we will work with the above normalized $\mathbf{h}_{ub}$. Consequently, we also need to normalize $\mathbf{h}^{(2)}_{ub}$. Thus, henceforth, $\mathbf{h}^{(2)}_{ub}\coloneqq \mathbf{h}^{(2)}_{ub}\frac{|\Psi(\eta\xi,\pi)|}{|c(\xi)|}$, where $\mathbf{h}^{(2)}_{ub}$ on the right is as in the previous section. 
\end{remark}

\begin{remark}
    Note that, when $\Sigma_{\infty}=\Sigma_{b}$ resp. $\Sigma_{\infty}=\Sigma_{ub}$, we may and will take $\eta=-1$ resp. $\eta=1$. This is possible since we are assuming $|\kappa|\equiv |\Sigma_{b}|\ (\text{mod}\ 2)$. In particular, when $\Sigma_{\infty}=\Sigma_{b}$ resp. $\Sigma_{\infty}=\Sigma_{ub}$, we may assume that, $\mathbf{h}_{ub}=\mathbf{h}$ resp. $\mathbf{h}_{ub}=\overline{\mathbf{h}}$. We will use this observation in the next subsection.
\end{remark}

\begin{remark}
   One can show that $\Psi(\eta\xi,\pi)$ is a totally real algebraic number (cf. \cite{HI13}, Defn. 4.1). Hence, the Fourier coefficients of $\mathbf{h}_{ub}$ are all totally real algebraic numbers, at least when $\Sigma_{b}$ or $\Sigma_{ub}$ is empty. We will use this fact later. 
\end{remark}

Let $\kappa_{0}=\max_{v\in\Sigma_{\infty}} \{\kappa_{v}\}$. For $\sigma\in \Aut(\mathbb{C})$, let $\pi^{\sigma}$ be the automorphic representation of $\PGL_{2}(\mathbb{A})$ as defined in \cite{RT11}, Theorem 3.10; in particular, weight of $\pi^{\sigma}$ equals $(2\kappa_{\sigma^{-1}\circ v})_{v\in \Sigma_{\infty}}$. Furthermore, if $r\in \{0,1\}^{n}$, let $u(r,\pi^{\sigma})$ denote Shimura's $u$-invariant (cf. \cite{Shi78}, Thm. 4.3; here it is understood that $u(r,\pi^{\sigma})\coloneqq u(r,f^{\sigma})$ where $f$ is the modular form corresponding to $\pi$).  
\begin{lemma}\label{ratmf}
For any $\sigma\in \Aut(\mathbb{C})$, we have
\begin{equation}\label{piLvaluerat}
    \Bigg(\frac{L(\frac{1}{2},\pi \otimes\chi_{\eta \xi})}{\pi ^{|\kappa_{0}-\kappa|}|\N_{F/\mathbb{Q}}\eta\xi|^{1/2}|u(\epsilon,\pi)|}\Bigg)^{\sigma}\mspace{-10mu}
    =\mspace{-6.0mu}\frac{L(\frac{1}{2},\pi^{\sigma} \otimes\chi_{\eta \xi})}{\pi ^{|\kappa_{0}-\kappa|}|\N_{F/\mathbb{Q}}\eta\xi|^{1/2}|u(\epsilon,\pi^{\sigma})|},
\end{equation}
where $\epsilon\in \{0,1\}^{n}$ is prescribed by the condition $\chi_{\eta\xi}(a)=\sgn{[a^{\epsilon}\N(a)^{\kappa_{0}}]}$ and $\N_{F/\mathbb{Q}}$ denotes the norm map.
\end{lemma}

\begin{proof}
 It is well known that $L_{\infty}(\frac{1}{2},\pi \otimes\chi_{\eta \xi})=\prod_{v\in\Sigma_{\infty}}(2\pi)^{-\kappa_{v}}\Gamma(\kappa_{v})$ and that $L(\frac{1}{2},\pi \otimes\chi_{\eta \xi})$ is a real number. The result now follows from the rationality of central values of $L$-functions associated to automorphic representations of Hilbert modular forms given in Theorem 4.3 of \cite{Shi78} and a result on Gauss sums (\cite{Shi87}, Lemma 9.3).  
\end{proof}

We thus have the following result on the rationality of $\frac{L(1,\pi,\Ad)}{\langle\mathbf{h}_{ub},\mathbf{h}_{ub}\rangle}$.
\begin{prop}\label{secondfac}
For any $\sigma\in \Aut(\mathbb{C})$,
\begin{equation}
    \Bigg(\frac{L(1,\pi,\Ad)}{\pi ^{|\kappa_{0}-\kappa+1|}\langle\mathbf{h}_{ub},\mathbf{h}_{ub}\rangle |u(\epsilon,\pi)|}\Bigg)^{\sigma}=\frac{L(1,\pi^{\sigma},\Ad)}{\pi ^{|\kappa_{0}-\kappa+1|}\langle(\mathbf{h}_{ub})^{\sigma},(\mathbf{h}_{ub})^{\sigma}\rangle |u(\epsilon,\pi^{\sigma})|},
\end{equation}
where $\epsilon$ is as in lemma \ref{ratmf}. 
\end{prop}
\begin{proof}
Follows from the previous two lemmas together with the well known fact that \begin{equation}\label{wellknown}
    D_{F}^{1/2}\xi_{F}(2)\in \pi^{n}\mathbb{Q}.
\end{equation}
\end{proof}

\begin{remark}
    As already mentioned in a previous Remark, we may and will take $\eta=-1$ resp. $\eta=1$ when $\Sigma_{\infty}=\Sigma_{b}$ resp. $\Sigma_{\infty}=\Sigma_{ub}$. Moreover, in these cases we have, respectively, $\epsilon=(\kappa_{0}+1,\ldots,\kappa_{0}+1)\ (\text{mod}\ 2)$ and $\epsilon=(\kappa_{0},\ldots,\kappa_{0})\ (\text{mod}\ 2)$. 
\end{remark}

Next, we will work out the rationality of the global period $|\langle \mathbf{g}^{(\delta)}_{ub}, (Y_{+}^{(r)}F_{\mathbf{h}_{ub}})|_{\SL_{2}(\mathbb{A})}\rangle|^{2}$ in the two extreme cases, the \textit{purely balanced} and the \textit{purely unbalanced}, and as a consequence derive the rationality of the central value $L(\frac{1}{2},\Sym^{2}\tau\times \pi)$ in these cases. We will make essential use of results from Shimura's papers \cite{Shi78} and \cite{Shi87}. 

\begin{subsection}{The Purely Balanced Case}
We refer to the case when $\Sigma_{\infty}=\Sigma_{b}$ as purely balanced. In this case, as already noted in the second remark after Lemma \ref{cor3.3}, we may take $\mathbf{h}_{ub}=\mathbf{h}$. Consequently, $Y_{+}^{(r)}F_{\mathbf{h}_{ub}}=Y_{+}^{(r)}F_{\mathbf{h}}$. Moreover, since $\Sigma_{ub}=\emptyset$, we may denote $\mathbf{g}^{(\delta)}_{ub}$ by $\mathbf{g}^{(\delta)}$.

Now, $(Y_{+}^{(r)}F_{\mathbf{h}})\big|_{\SL_{2}(\mathbb{A})}$ is a nearly holomorphic Hilbert modular form of integral weight by Lemma \ref{nearhollemma}, since $\X_{-}^{(r)}(Y_{+}^{(r)}F_{\mathbf{h}})=0.$
Denote its holomorphic projection by $\mathbf{g}_{0}$. Then by (\cite{Shi87}, Prop. 9.4(i)),
\begin{equation*}
    \langle \mathbf{g}^{(\delta)}, (Y_{+}^{(r)}F_{\mathbf{h}})\big|_{\SL_{2}(\mathbb{A})}\rangle=\langle \mathbf{g}^{(\delta)},\mathbf{g}_{0}\rangle.
\end{equation*}

We claim that the Fourier coefficients of $\mathbf{g}^{(\delta)}$ and $\mathbf{g}_{0}$ are totally real. That the Fourier coefficients of $\mathbf{g}^{(\delta)}$ are totally real follows from the fact that the central character of $\tau$ is trivial. For $\mathbf{g}_{0}$, we begin by noting that, for any $\sigma\in \Aut(\mathbb{C})$, 
$$((Y_{+}^{(r)}F_{\mathbf{h}})\big|_{\SL_{2}(\mathbb{A})})^{\sigma}=(Y_{+}^{(r),\sigma}F_{\mathbf{h}^{\sigma}})\big|_{\SL_{2}(\mathbb{A})},$$ 
where $Y_{+}^{(r),\sigma}$ is the adelization of the differential operator $\prod_{v\in\Sigma_{\infty}}\Delta_{v}^{r_{\sigma\circ v}}$ introduced in the Introduction (when $F=\mathbb{Q}$, this has been proved in \cite{Xue19}, Lemma 2.3). Denoting the holomorphic projection operator by $pr_{hol}$, it follows that
$$pr_{hol}((Y_{+}^{(r),\sigma}F_{\mathbf{h}^{\sigma}})\big|_{\SL_{2}(\mathbb{A})})=pr_{hol}(((Y_{+}^{(r)}F_{\mathbf{h}})\big|_{\SL_{2}(\mathbb{A})})^{\sigma})=pr_{hol}((Y_{+}^{(r)}F_{\mathbf{h}})\big|_{\SL_{2}(\mathbb{A})})^{\sigma}=\mathbf{g}_{0}^{\sigma},$$
where the second equality follows from (\cite{Shi87}, Prop. 9.4(ii)). Recall that by our normalization, the Fourier coefficients of $\mathbf{h}$ are totally real (see the third Remark after Lemma 4.1); thus by the above equality the Fourier coefficients of $\mathbf{g}_{0}$ are also totally real. 
Therefore,
\begin{equation}\label{simabsval}
    |\langle \mathbf{g}^{(\delta)}, (Y_{+}^{(r)}F_{\mathbf{h}})\big|_{\SL_{2}(\mathbb{A})}\rangle|^{2}=\langle \mathbf{g}^{(\delta)},\mathbf{g}_{0}\rangle^{2}.
\end{equation}

\begin{prop}\label{innprodrat}
For any $\sigma\in \Aut(\mathbb{C})$,
\begin{equation*}
    \sigma\Bigg(\frac{|\langle \mathbf{g}^{(\delta)}, (Y_{+}^{(r)}F_{\mathbf{h}})\big|_{\SL_{2}(\mathbb{A})}\rangle|^{2}}{\langle \mathbf{g}
    ,\mathbf{g}
    \rangle^{2}}\Bigg)=\frac{|\langle (\mathbf{g}^{(\delta)})^{\sigma}, (Y_{+}^{(r),\sigma}F_{\mathbf{h}^{\sigma}})\big|_{\SL_{2}(\mathbb{A})}\rangle|^{2}}{\langle \mathbf{g}^{\sigma},\mathbf{g}^{\sigma}\rangle^{2}}.
\end{equation*}

\end{prop}
\begin{proof}
Let $\mathbf{g}$ denote the newform $\otimes_{v}g_{v}$ with $g_{v}$'s as defined in $\S\ref{loccomp}$. Then we can write $\langle \mathbf{g}^{(\delta)},\mathbf{g}_{0}\rangle$ as $\langle \mathbf{g},\mathbf{g}_{0}'\rangle$ for some $\mathbf{g}_{0}'$ since we are working with unitary representations. Now, by (\cite{Shi78}, Prop. 4.15), for any  $\sigma\in\Aut(\mathbb{C})$,
\begin{equation}
    \sigma\Bigg(\frac{\langle \mathbf{g},\mathbf{g}_{0}'\rangle}{\langle \mathbf{g}
    ,\mathbf{g}
    \rangle}\Bigg)=\frac{\langle \mathbf{g}^{\sigma},\mathbf{g}_{0}'^{\sigma}\rangle}{\langle \mathbf{g}^{\sigma},\mathbf{g}^{\sigma}\rangle}.
\end{equation}
The proposition now follows from (\ref{simabsval}).
\end{proof}

\begin{theorem}\label{ratb}
Let $\pi$ and $\tau$ be irreducible cuspidal automorphic representations of $\GL_{2}(\mathbb{A})$ of weights $2\kappa$ and $\kappa'+1$ respectively. Assume that both $\pi$ and $\tau$ are of full level. Then for any $\sigma\in \Aut(\mathbb{C})$,
\begin{equation*}
    \sigma\Bigg(\frac{L(\frac{1}{2},\Sym^{2}\tau\times\pi)}{\pi ^{|\kappa_{0}-\kappa+2|}\sqrt{D_{F}}\langle \mathbf{g},\mathbf{g}\rangle^{2}|u(\epsilon,\pi)|}\Bigg)=\frac{L(\frac{1}{2},\Sym^{2}\tau^{\sigma}\times\pi^{\sigma})}{\pi ^{|\kappa_{0}-\kappa+2|}\sqrt{D_{F}}\langle \mathbf{g}^{\sigma},\mathbf{g}^{\sigma}\rangle^{2}|u(\epsilon,\pi^{\sigma})|},
\end{equation*}
where $\mathbf{g}$ is the newform generating $\tau$ and $\epsilon=(\kappa_{0}+1,\ldots,\kappa_{0}+1)\ (\text{mod}\ 2)$. In particular,
\begin{equation*}
    L(\frac{1}{2},\Sym^{2}\tau\times\pi)\sim_{\mathbb{Q}(\pi,\tau)}\pi^{3n}\sqrt{D_{F}}\langle \mathbf{g},\mathbf{g}\rangle^{2}\langle \mathbf{h},\mathbf{h}\rangle,
\end{equation*}
where $\mathbb{Q}(\pi,\tau)$ is the number field generated by the Fourier coefficients of the Hilbert modular forms associated to $\pi$ and $\tau$, and $\sim_{\mathbb{Q}(\pi,\tau)}$ means up to an element of $\mathbb{Q}(\pi,\tau)$. 
\end{theorem}

\begin{proof}
The first and last assertions follow by combining Prop. \ref{innprodrat}, Prop. \ref{secondfac} (with $\epsilon=(\kappa_{0}+1,\ldots,\kappa_{0}+1)\ (\text{mod}\ 2)$) and fact \ref{wellknown} with Thm. \ref{mainthm}. The second assertion follows from the first by noting that eqn. (\ref{uninnprodrat}) below implies that 
$$\pi^{|\kappa_{0}-\kappa-1|}|u(\epsilon,\pi)|\sim_{\mathbb{Q}(\pi)^{\times}}\langle \mathbf{h}^{(2)},\mathbf{h}^{(2)}\rangle.$$
\end{proof}
\end{subsection}

\begin{subsection}{The Purely Unbalanced Case}
We refer to the case when $\Sigma_{\infty}=\Sigma_{ub}$ as purely unbalanced. In this case, as already noted in the second Remark after Lemma 4.1, we may assume $\mathbf{h}_{ub}=\overline{\mathbf{h}}$. Consequently, $Y_{+}^{(r)}F_{\mathbf{h}_{ub}}=Y_{+}^{(r)}F_{\overline{\mathbf{h}}}$. Moreover, we may denote $\mathbf{g}^{(\delta)}_{ub}$ by $\overline{\mathbf{g}}^{(\delta)}$. This time we will use the theory of nearly holomorphic Hilbert modular forms of \textit{half-integral} weight to work out the rationality of the central $L$-value. 

Note that, $\mathbf{h}^{(2)}_{ub}=\overline{\mathbf{h}}^{(2)}$ is anti-holomorphic and, by our choice of local components (\S\ref{loccomp}), so is $\overline{\mathbf{g}}^{(\delta)}$. Thus their complex conjugates are holomorphic vectors. Write $F_{\overline{h}}$ as $\overline{\mathbf{h}}^{(2)}\otimes \Theta_{Y_{+}^{(r)}\boldsymbol{\phi}}$. Then, it is not hard to see that,
\begin{equation*}
    \langle (Y_{+}^{(r)}F_{\overline{\mathbf{h}}})\big|_{\SL_{2}(\mathbb{A})},\overline{\mathbf{g}}^{(\delta)}\rangle=\langle \overline{\mathbf{h}}^{(2)}\otimes \Theta_{Y_{+}^{(r)}\boldsymbol{\phi}}\big|_{\SL_{2}(\mathbb{A})}, \overline{\mathbf{g}}^{(\delta)}\rangle=\langle  \mathbf{g}^{(\delta)}\otimes \Theta_{\Y_{+}^{(r)}\boldsymbol{\phi}}\big|_{\SL_{2}(\mathbb{A})}, \mathbf{h}^{(2)}\rangle,
\end{equation*}
where $\mathbf{g}^{(\delta)}\otimes Y_{+}^{(r)}\boldsymbol{\phi}\big|_{\SL_{2}(\mathbb{A})}$ is now \textit{nearly holomorphic of half-integral weight} by Lemma \ref{nearhollemma}. Let $\mathbf{h}_{0}$ be its holomorphic projection. Then, by (\cite{Shi87}, Prop. 9.4),
\begin{equation*}
    \langle  \mathbf{g}^{(\delta)}\otimes \Theta_{Y_{+}^{(r)}\boldsymbol{\phi}}\big|_{\SL_{2}(\mathbb{A})}, \mathbf{h}^{(2)}\rangle=\langle \mathbf{h}_{0}, \mathbf{h}^{(2)}\rangle.
\end{equation*}
By a reasoning similar to that in the purely balanced case, the Fourier coefficients of $\mathbf{g}^{(\delta)}, \mathbf{h}^{(2)}$ and $\mathbf{h}_{0}$ are all totally real. Therefore,
\begin{equation}\label{simabsvalub}
    |\langle \overline{\mathbf{h}}^{(2)}\otimes \Theta_{Y_{+}^{(r)}\boldsymbol{\phi}}\big|_{\SL_{2}(\mathbb{A})}, \overline{\mathbf{g}}^{(\delta)}\rangle|^{2}=|\langle  \mathbf{g}^{(\delta)}\otimes \Theta_{Y_{+}^{(r)}\boldsymbol{\phi}}\big|_{\SL_{2}(\mathbb{A})}, \mathbf{h}^{(2)}\rangle|^{2}=\langle \mathbf{h}_{0}, \mathbf{h}^{(2)}\rangle^{2}.
\end{equation}
\begin{prop}\label{innprodratub}
For all $\sigma\in \Aut(\mathbb{C})$,
\begin{equation*}
    \sigma\Bigg(\frac{|\langle (Y_{+}^{(r)}F_{\overline{\mathbf{h}}})\big|_{\SL_{2}(\mathbb{A})}, \overline{\mathbf{g}}^{(\delta)}\rangle|^{2}}{\pi^{2|\kappa_{0}-\kappa-1|}u(\epsilon,\pi)^{2}}\Bigg)=\frac{|\langle (Y_{+}^{(r),\sigma}F_{\overline{\mathbf{h}}^{\sigma}})\big|_{\SL_{2}(\mathbb{A})}, (\overline{\mathbf{g}}^{(\delta)})^{\sigma}\rangle|^{2}}{\pi^{2|\kappa_{0}-\kappa-1|}u(\epsilon,\pi^{\sigma})^{2}},
\end{equation*}
where $\epsilon=(\kappa_{0}+1,\ldots,\kappa_{0}+1)\ (\text{mod}\ 2)$ and $Y_{+}^{(r),\sigma}$ is the adelization of the differential operator $\prod_{v\in\Sigma_{\infty}}\Delta_{v}^{r_{\sigma\circ v}}$ introduced in the Introduction.
\end{prop}
\begin{proof}
By (\cite{Shi87}, Thm. 10.5), for all $\sigma\in\Aut(\mathbb{C})$
\begin{equation}\label{uninnprodrat}
    \sigma\Bigg(\frac{\langle \mathbf{h}_{0}, \mathbf{h}^{(2)}\rangle}{\pi^{|\kappa_{0}-\kappa-1|}|u(\epsilon,\pi)|}\Bigg)=\frac{\langle \mathbf{h}_{0}^{\sigma}, (\mathbf{h}^{(2)})^{\sigma}\rangle}{\pi^{|\kappa_{0}-\kappa-1|}|u(\epsilon,\pi^{\sigma})|},
\end{equation}
with $\epsilon$ as above.
The proposition now follows from (\ref{simabsvalub}). 
\end{proof}

\begin{theorem}\label{ratub}
Let $\pi$ and $\tau$ be irreducible cuspidal automorphic representations of $\GL_{2}(\mathbb{A})$ of weights $2\kappa$ and $\kappa'+1$ respectively. Assume that both $\pi$ and $\tau$ are of full level. Then for any $\sigma\in \Aut(\mathbb{C})$,
\begin{equation*}
    \sigma\Bigg(\frac{L(\frac{1}{2},\Sym^{2}\tau\times\pi)}{\pi ^{|\kappa_{0}-\kappa+1|}\sqrt{D_{F}}\langle \mathbf{f},\mathbf{f}\rangle |u(\epsilon,\pi)|}\Bigg)=\frac{L(\frac{1}{2},\Sym^{2}\tau^{\sigma}\times\pi^{\sigma})}{\pi ^{|\kappa_{0}-\kappa+1|}\sqrt{D_{F}}\langle \mathbf{f}^{\sigma},\mathbf{f}^{\sigma}\rangle |u(\epsilon,\pi^{\sigma})|},
\end{equation*}
where $\mathbf{f}$ is the newform generating $\pi$ and $\epsilon=(\kappa_{0}+1,\ldots,\kappa_{0}+1)\ (\text{mod}\ 2)$. In particular,
\begin{equation*}
    L(\frac{1}{2},\Sym^{2}\tau\times\pi)\sim_{\mathbb{Q}(\pi)}\pi^{2n}\sqrt{D_{F}}\langle \mathbf{f},\mathbf{f}\rangle\langle \mathbf{h},\mathbf{h}\rangle,
\end{equation*}
where $\mathbb{Q}(\pi)$ is the number field generated by the Fourier coefficients of the Hilbert modular form corresponding to $\pi$ and $\sim_{\mathbb{Q}(\pi)}$ means up to an element of $\mathbb{Q}(\pi)$. 
\end{theorem}
\begin{proof}
It is well-known that (cf. \cite{Shi87}, Thm 10.2(II)) 
\begin{equation}\label{piinnprod}
    \sigma\Bigg(\frac{u(\epsilon,\pi)u(1+\epsilon,\pi)}{(2\pi i)^{|2\kappa-2\kappa_{0}+1|}\langle\mathbf{f},\mathbf{f}\rangle}\Bigg)=\frac{u(\epsilon,\pi^{\sigma})u(1+\epsilon,\pi^{\sigma})}{(2\pi i)^{|2\kappa-2\kappa_{0}+1|}\langle\mathbf{f}^{\sigma},\mathbf{f}^{\sigma}\rangle}.
\end{equation}
The first and last assertions follow by combining Prop. \ref{innprodratub}, Prop. \ref{secondfac} (with $\epsilon$ in the Prop. equal to $(\kappa_{0},\ldots,\kappa_{0})\ (\text{mod}\ 2)$), Eqn. (\ref{piinnprod}) and fact \ref{wellknown} with Thm. \ref{mainthm}. The second assertion follows from the first by noting that (\ref{uninnprodrat}) implies
$$\pi^{|\kappa_{0}-\kappa-1|}|u(\epsilon,\pi)|\sim_{\mathbb{Q}(\pi)}\langle \mathbf{h}^{(2)},\mathbf{h}^{(2)}\rangle$$
and that $\langle \mathbf{h}^{(2)},\mathbf{h}^{(2)}\rangle=\langle \mathbf{h},\mathbf{h}\rangle$
\end{proof}

\begin{remark}
Comparing Thm. \ref{ratb} and Thm. \ref{ratub}, we see that, the period in the purely balanced case depends on both $\pi$ and $\tau$, whereas in the purely unbalanced case depends only on $\pi$.
\end{remark}
\end{subsection}
\end{section}

\begin{section}{A Conjecture For The General Case}\label{generalcase}
In this section we present a conjecture related to the rationality of $L(\frac{1}{2},\Sym^{2}\tau\times \pi)$ for the general case, that is, when the number of balanced and unbalanced places are arbitrary. In fact, we will state the cconjecture not only for the central value but for all critical values of $L(s,\Sym^{2}\tau\times \pi)$.

A critical point for $L(s,\Sym^{2}\tau\times \pi)$ is a half-integer $m+\frac{1}{2}$ such that neither $L_{\infty}(s,\Sym^{2}\tau\times \pi)$ nor $L_{\infty}(1-s,\Sym^{2}\tau\times \pi)$ has a pole at $s=m+\frac{1}{2}$. Recall that \begin{equation*}
    r_{v}=\begin{cases}
    \kappa_{v}-\kappa'_{v}-1 & \text{if}\ v\in \Sigma_{ub}\\
    \kappa'_{v}-\kappa_{v} & \text{if}\ v\in \Sigma_{b}.
    \end{cases}
\end{equation*}
Put 
$$t^{0}=\min_{v\in \Sigma_{ub}}\{r_{v}\},\ \ \ a^{0}=\min_{v\in \Sigma_{b},2\kappa_{v}>\kappa_{v}'}\{r_{j}\},\ \ \ b^{0}=\min_{v\in \Sigma_{b},\kappa'\geq 2\kappa}\{\kappa_{v}-1\}_{v\in \Sigma_{b},\kappa'\geq 2\kappa}.$$
\begin{prop}
 The set of critical points for $L(s,\Sym^{2}\tau\times \pi)$ is
 \begin{equation*}
     \Big\{m+\frac{1}{2}\in \mathbb{Z}+\frac{1}{2}\mid -\min\{a^{0},b^{0},t^{0}\}\leq m \leq\min\{a^{0},b^{0},t^{0}\}\Big\}
 \end{equation*}
In particular, $s=\frac{1}{2}$ is a critical point.
\end{prop}
\begin{proof}
Follows from the definition of a critical point and Eq. (\ref{linfty}) after a tedious but elementary computation. 
\end{proof}
Deligne's conjecture on special values of $L$-functions predicts the transcendental part of the value of a (motivic) L-function at critical points. It formulates the transcendental part of such values in terms of certain complex numbers known as Deligne's periods which are obtained by comparing the Betti and the deRham cohomological realizations of the underlying motive. The following conjecture is suggested by explicitly computing Deligne's period in our case using results from (\cite{Yosh94},\S2) and \cite{Yosh95}. 
\begin{conjecture}\label{mcp}
If $m+\frac{1}{2}\in \mathbb{Z}+\frac{1}{2}$ is a critical point of the automorphic $L$-function $L(s,\Sym^{2}\tau\times \pi)$, then, for all $\sigma\in \Aut(\mathbb{C})$,
\begin{equation}
\resizebox{0.99\hsize}{!}
   {$\Big(\frac{L(m+\frac{1}{2},\Sym^{2}\tau\times \pi)}{\pi ^{|\kappa_{0}-\kappa+1|+|\Sigma_{b}|}\sqrt{D_{F}}|u(\epsilon,\pi)|\nu^{\Sigma_{ub}}(\pi)\nu^{\Sigma_{b}}(\tau)^{2}}\Big)^{\sigma}= \frac{L(m+\frac{1}{2},\Sym^{2}\tau^{\sigma}\times \pi^{\sigma})}{\pi ^{|\kappa_{0}-\kappa+1|+|\Sigma_{b}|}\sqrt{D_{F}}|u(\epsilon,\pi^{\sigma})|\nu^{\Sigma_{ub}}(\pi^{\sigma})\nu^{\Sigma_{b}}(\tau^{\sigma})^{2}}$},
\end{equation}
where $\epsilon=(\epsilon_{v})$ with $\epsilon_{v}=m+\kappa_{0}+1\ (\text{mod}\  2)\ \  \text{for all $v$}$, and $u(\cdot\  ,\ \cdot )$ and $\nu^{\cdot}(\ \cdot\ )$ are respectively, Shimura's $u$-invariant and Harris' period defined in (\cite{Har90}, \S1). 
\end{conjecture}

\begin{remark}
    The above conjecture is compatible with Theorems \ref{ratb} and \ref{ratub} since $\nu^{\Sigma_{b}}(\tau)\sim_{\mathbb{Q}(\tau)^{\times}}\langle \mathbf{g},\mathbf{g}\rangle$ when $\Sigma_{b}=\Sigma_{\infty}$, and $\nu^{\Sigma_{ub}}(\pi)\sim_{\mathbb{Q}(\pi)^{\times}}\langle \mathbf{f},\mathbf{f}\rangle$ when $\Sigma_{ub}=\Sigma_{\infty}$.
\end{remark}

When $m=0$, Prop. \ref{secondfac} together with properties of Shimura's invariants (cf. \cite{Yosh95}) suggest that the above conjecture is equivalent to 
\begin{conjecture}
For every $S\subset\Sigma_{\infty}, r\in\{0,1\}^{S}$ and $\sigma\in\Aut(\mathbb{C})$ there exists a complex number $P(\pi^{\sigma},S,r)$ such that, for any $\sigma\in \Aut(\mathbb{C})$,
\begin{equation*}
    \Bigg|\frac{\langle\mathbf{g}^{(\delta)}_{ub},(Y_{+}^{(r)}F_{\mathbf{h}_{ub}})\big|_{\SL_{2}(\mathbb{A})}\rangle}{P(\pi,\Sigma_{ub},\epsilon)\nu^{\Sigma_{b}}(\tau)}\Bigg|^{\sigma}=\Bigg|\frac{\langle(\mathbf{g}^{(\delta)}_{ub})^{\sigma},(Y_{+}^{(r),\sigma}F_{\mathbf{h}^{\sigma},ub})\big|_{\SL_{2}(\mathbb{A})}\rangle}{{P(\pi^{\sigma},\Sigma_{ub}^{\sigma},\epsilon^{\sigma})}\nu^{\Sigma_{b}^{\sigma}}(\tau^{\sigma})}\Bigg|,
\end{equation*}
where $\epsilon=(\epsilon_{v})$ with $\epsilon_{v}=m+\kappa_{0}+1\ (\text{mod}\  2)\ \  \text{for all $v\in \Sigma_{ub}$}$, $\epsilon^{\sigma}=(\epsilon_{\sigma\circ v})$ and $S^{\sigma}=\{\sigma\circ v| v\in S\}$.
The complex number $P(\pi^{\sigma},S,r)$ should be thought of as a cohomological interpretation of Shimura's $P$-invariant defined in \cite{Yosh95}.
\end{conjecture}

In the remainder of this paper we will prove Proposition \ref{localperiod}.
\end{section}

\begin{section}{Computation at the Archimedean Places}\label{archplaces}
\begin{subsection}{The Unbalanced Case}

Let $v\in\Sigma_{ub}$. For simplicity we will suppress all subscripts $v$ throughout this section. Recall (ref. \S\ref{loccomp}) the following notation at $v$. 
\begin{itemize}
    \item $\tau$ is a discrete series representation of $\PGL_{2}(\mathbb{R})$ of mimimal weight $\pm(\kappa'+1)$ and $g$ (resp. $\overline{g}$)  is a lowest (resp. highest weight vector) in $\tau$. 
    \item $\sigma$ is a discrete series representation of $\widetilde{\SL_{2}(\mathbb{R})} $ of highest weight $-(\kappa+\frac{1}{2})$ and $\overline{h}\in \sigma$ is a highest weight vector. Furthermore, let $\overline{\sigma}$ denote the contragradient of $\sigma$ and $h\in \overline{\sigma}$ be the holomorphic vector corresponding to $h$.
    \item $\psi$ is the additive character of $\mathbb{R}$ given by $x\mapsto e^{2\pi\imag x}$. Furthermore, $\omega_{\psi}$ is the Schr\"{o}dinger-Weil representation of $\SL_{2}(\mathbb{R})$ and $\phi=e^{-2\pi\imag x^{2}}\in \omega_{\psi}$ is a lowest weight vector of weight $1/2$.
    \item We have $\kappa> \kappa'$ and $r=\kappa- \kappa'-1.$
\end{itemize}
The goal of this subsection is to compute
$$\Pe^{\sharp}(\overline{g},\overline{h},Y_{+}^{r}\phi)=\Bigg(\frac{\xi_{F}(2)L(1/2,\Sym^{2}\tau\times \pi)}{L(1,\pi,\Ad)L(1,\tau,\Ad)}\Bigg)^{-1}\frac{\Pe(\overline{g},\overline{h},Y_{+}^{r}\phi)}{ \langle \overline{g},\overline{g}\rangle \langle \overline{h},\overline{h}\rangle \langle \phi,\phi\rangle},$$
where  $$\Pe( \overline{g}, \overline{h},Y_{+}^{r}\phi)=\int_{\SL_{2}(\mathbb{R})}\langle \tau(X)\overline{g},\overline{g}\rangle\overline{\langle \sigma(X)\overline{h},\overline{h}\rangle \langle \omega_{\psi}(X)Y_{+}^{r}\phi,Y_{+}^{r}\phi\rangle}  dX.$$
\begin{prop}\label{archk>k'}
We have,
$$\Pe^{\sharp}(\overline{g}, \overline{h},Y_{+}^{r}\phi)=\frac{2^{2r_v+1}}{r_v!}\pi^{r_v}.$$
\end{prop}
Firstly, note that we can rewrite
\begin{align*}
    \Pe( \overline{g}, \overline{h},Y_{+}^{r}\phi)
    &=\int_{\SL_{2}(\mathbb{R})}\overline{\langle \tau(X)g,g\rangle}\langle \overline{\sigma}(X)h,h\rangle \overline{\langle \omega_{\psi}(X)Y_{+}^{r}\phi,Y_{+}^{r}\phi\rangle}  dX,
\end{align*}

Now, the idea of the proof of the above proposition is to realize the holomorphic discrete series $\overline{\sigma}$ of weight $\kappa+\frac{1}{2}$ as a subrepresentation of $\tau\otimes\omega_{\psi}|_{\widetilde{\SL_{2}(\mathbb{R}})}$. 

We make use of explicit models for the representations $ \tau$ and $\omega_{\psi}$. For the full description of these models we refer the readers to (\cite{Bump}, Prop. 2.5.4) and (\cite{BS98}, Prop. 3.2.3), respectively.
Let the vector space underlying the model for $\omega_{\psi}$ be
$\bigoplus_{j\in \mathbb{N}_{0}}\mathbb{C}v_{j}$. 
Then the action of $\Z$ and $\X_{-}$ is given by 
\begin{equation}\label{omega}
  \Z v_{j}=(j+\frac{1}{2})v_{j},\ \ \ \ 
    \X_{-}v_{j}=\pi j(j-1)v_{j-2}.
\end{equation}
The vector space underlying the model for $\tau$ is  
$DS_{\pm(\kappa'+1)}=\bigoplus_{i\in \pm{2\mathbb{N}_{0}}}\mathbb{C}u_{i}$
and the action of $\Z$ and $\X_{-}$ is given by 
\begin{equation}\label{tau}
   \Z u_{i}=(i+\kappa'+1)u_{i},\ \ \ \ \X_{-}u_{i}=-\frac{i}{2}u_{i-2}. 
\end{equation}

Tensoring the representations $\tau$ and $\omega_ {\psi}|_{\widetilde{\SL_{2}(\mathbb{R}})}$, we get a (genuine) representation of $\widetilde{\SL_{2}(\mathbb{R}})$ whose underlying vector space is 
$$D_{\kappa'+\frac{3}{2}}=\bigoplus_{l\in 2\mathbb{N}_{0},\ k\in \mathbb{N}_{0}}\mathbb{C}u_{l}\otimes v_{k},$$
and the action of $\Z$ and $\X_{-}$ can be calculated using (\ref{omega}) and (\ref{tau});
$$\Z(u_{i}\otimes v_{j})=(i+j+\kappa'+\frac{3}{2})u_{i}\otimes v_{j},\ \ \ \ \ \ \ \
\X_{-}(u_{i}\otimes v_{j})=-\frac{i}{2}u_{i-2}\otimes v_{j}+\pi j(j-1)u_{i}\otimes v_{j-2}.$$

There is an inner product on $D_{\kappa'+\frac{3}{2}}$ such that $(u_{i}\otimes v_{j})$'s form an orthogonal basis. Denote this inner product by $\langle -,-\rangle$ and write $\langle v, v\rangle=||v||^{2}$. Then by (\cite{BS98}, p. 46-47) it follows
$$||u_{i}\otimes v_{j+1}||^{2}=2\pi (j+1)||u_{i}\otimes v_{j}||^{2},\ \ \ \ \\ (i+2\kappa')||u_{i}\otimes v_{j}||^{2}=i||u_{i-2}\otimes v_{j}||^{2}.$$
We may normalize the inner product so that $||u_{0}\otimes v_{t}||=1$. Then for any $2\leq i\leq r $, $i$ even, we have
\begin{equation}\label{ipformula}
  ||u_{i}\otimes v_{r-i}||^{2}=(2\pi)^{-i}\prod_{\substack{0\leq l\leq i-2 \\ l\in 2\mathbb{N}_{0}}}\frac{(i-l)}{(i+2\kappa'-l)(r-i+l+1)(r-i+l+2)}.  
\end{equation}

The space
$$D_{\kappa'+\frac{3}{2}}(r)=\bigoplus_{i+j=r,i\in 2\mathbb{N}_{0},j\in \mathbb{N}_{0}}\mathbb{C}u_{i}\otimes v_{j}$$
is the largest subspace on which $\Z$ acts by $r+\kappa'+\frac{3}{2}=\kappa+\frac{1}{2}$.
\begin{lemma}\label{v}
There is a unique (upto a scalar) vector $v_{t}^{\hol}$ in $D_{\kappa'+\frac{3}{2}}(r)$ with the property that $\X_{-}v_{t}^{\hol}=0.$ It is given by 
\begin{equation*}
   v_{r}^{\hol}=\sum_{\substack{0\leq i\leq r \\ i\in 2\mathbb{N}_{0}}}c_{i}u_{i}\otimes v_{r-i};\ \ \\ \ c_{0}=1,\ \ \ \ c_{i}=(2\pi )^{i/2}\prod_{\substack{0\leq l\leq i-2 \\ l\in 2\mathbb{N}_{0}}}\frac{(r-l)(r-l-1)}{(i-l)},\ \ (i\geq 2). 
\end{equation*}
In particular,
\begin{equation}\label{||v||^{2}}
    ||v_{r}^{\hol}||^{2}=2^{-r}\frac{{2\kappa-2\choose r}}{{\kappa-1 \choose r}}
\end{equation}
\end{lemma}
\begin{proof}
Suppose that 
$$v_{r}^{\hol}=\sum_{\substack{0\leq i\leq r \\ i\in 2\mathbb{N}_{0}}}c_{i}u_{i}\otimes v_{r-i}$$
and $\X_{-}v_{r}^{\hol}=0.$ Then by the formula for action of $\X_{-}$, we conclude that for any $0\leq i\leq r-2$ and $i$ even, we have
$$-c_{i+2}\frac{(i+2)}{2}=\pi (r-i)(r-i-1)c_{i}.$$
Let $c_{0}=1$. Then we may recursively solve for $c_{i}$'s. Furthermore,
\begin{align*}
    ||v_{r}^{\hol}||^{2}&=\sum_{\substack{0\leq i\leq r\\ i\in 2\mathbb{N}_{0}}}|c_{i}|^{2}||u_{i}\otimes v_{r-i}||^{2}\\
    &=1+\sum_{\substack{2\leq i\leq r \\ i\in 2\mathbb{N}_{0}}}\prod_{\substack{0\leq l\leq i-2 \\ l\in 2\mathbb{N}_{0}}} \frac{(r-l)^{2}(r-l-1)^{2}}{(i-l)(i+2\kappa-2r-l-2)(r-i+l+1)(r-i+l+2)}
\end{align*}
The lemma now follows from Lemma \ref{Monsterid2} given in the appendix.
\end{proof}
Thus, we realize $\overline{\sigma}$ as a subrepresentation of $\tau\otimes\omega_{\psi}|_{\widetilde{\SL_{2}(\mathbb{R}})}$ generated by $v_{r}^{\hol}$. We may assume that the inner product on $\overline{\sigma}$ is given by the restriction of that of $\tau\otimes\omega_{\psi}|_{\widetilde{\SL_{2}(\mathbb{R}})}$. Since $\Pe_{v}^{\sharp}( \overline{g}, \overline{h},Y_{+}^{r}\mathbf{\phi})$ does not change if we replace $ h,\phi$ or $ g$ by a scalar multiple of them, we may assume that $ h=v_{r}^{\hol}$ and  $ g\otimes Y_{+}^{r}\phi=u_{0}\otimes v_{r}$.

Now we can prove the proposition. The orthogonal projection of $u_{0}\otimes v_{r}$ to the line generated by $v_{r}^{\hol}$ is $||v_{r}^{\hol}||^{-2}v_{r}^{\hol}$. It follows that
$$\Pe( \overline{g}, \overline{h},Y_{+}^{r}\mathbf{\phi})=\frac{1}{||v_{r}^{\hol}||^{4}}\int_{\SL_{2}(\mathbb{R})}|\langle\overline{\sigma}(X)v_{r}^{\hol},v_{r}^{\hol}\rangle|^{2}dX.$$
As $\overline{\sigma}$ is the discrete series representation of ${\widetilde{\SL_{2}(\mathbb{R}})}$ with lowest weight $\kappa+\frac{1}{2}$, it is known that (ref. \cite{Xue18}, Lemma 5.2),
$$|\langle \overline{\sigma}(\diag[e^{t},e^{-t}])v_{r}^{\hol},v_{r}^{\hol}\rangle|=||v_{r}^{\hol}||^{2}\times (\cosh t)^{-(\kappa+\frac{1}{2})},\ \\ \ \ \ r\geq 0.$$
Let $X=k_{1}\diag[e^{t},e^{-t}]k_{2}$ be the Cartan decomposition. Then $dX=2\pi \sinh 2tdtdk_{1}dk_{2}$,
where $dk_{1},dk_{2}$ are the measure on $\SO_{2}(\mathbb{R})$ so that its volume is one and $dt$ is the usual Lebesgue measure on $\mathbb{R}$. Therefore
\begin{equation}\label{someeqn}
    \Pe_{v}^{\sharp}( \overline{g}, \overline{h},Y_{+}^{r}\mathbf{\phi})
    =\Bigg(\frac{\xi(2)L(1/2,\Sym^{2}\tau\times \pi)}{L(1,\pi,\Ad)L(1,\tau,\Ad)}\Bigg)^{-1}\frac{\displaystyle\int_{0}^{\infty}(\cosh t)^{-(2\kappa+1)}2\pi \sinh 2tdt}{||v_{r}^{\hol}||^{2}}.
\end{equation}
By definition,
\begin{align*}
 \frac{\xi(2)L(1/2,\Sym^{2}\tau\times \pi)}{L(1,\pi,\Ad)L(1,\tau,\Ad)}&=\frac{\pi^{-1}\Gamma(1)\cdot 2^{2}(2\pi)^{-2\kappa}\Gamma(\kappa+\kappa')\Gamma(\kappa-\kappa')2(2\pi)^{-\kappa}\Gamma(\kappa)}{2(2\pi)^{-\kappa'-1}\Gamma(\kappa'+1)\pi^{-1}\Gamma(1)\cdot 2^{2}(2\pi)^{-2\kappa}\Gamma(2\kappa)\pi^{-1}\Gamma(1)}\\
&=(2\pi)^{1-r}\frac{t!}{2\kappa-1}\frac{{\kappa-1 \choose r}}{{2\kappa-2\choose r}}.
\end{align*}
Moreover (cf. \cite{GR00}, \S2.433, \#11),
$$
\int_{0}^{\infty}(\cosh t)^{-(2\kappa+1)}2\pi \sinh 2tdt=4\pi(2\kappa-1)^{-1}.$$
Proposition \ref{archk>k'} is now proved by plugging (\ref{||v||^{2}}) and the formulae above into (\ref{someeqn}).
\end{subsection}

\begin{subsection}{The Balanced Case}
\begin{prop}\label{archk'>k}
 Let $v\in \Sigma_{\infty}$ be such that $\kappa_{v}'\geq \kappa_{v}.$ Put $r_{v}=\kappa_{v}'- \kappa_{v}.$ Then 
 $$\Pe_{v}^{\sharp}( g_{v}, h_{v},Y_{+}^{r_{v}}\mathbf{\phi}_{v})=\frac{1}{r_{v}!}\pi^{r_{v}}.$$
\end{prop}

The proof of this proposition is similar to the one in the unbalanced case except that this time we will realize the \textit{integral} weight representation $\tau_{v}|_{\SL_{2}(\mathbb{R})}$ as a subrepresentation of $\sigma_{v}\otimes\omega_{\psi_{v}}|_{\SL_{2}(\mathbb{R})}$. It has been worked out in detail
in (\cite{Xue19},\S3.4) except that the rational part of $\Pe_{v}^{\sharp}$ is not explicitly computed. We will only supply here the information left out in \cite{Xue19}. 

In (\cite{Xue19}, Prop.3.4), it has been shown that
\begin{equation}\label{someeqn1}
    \Pe_{v}^{\sharp}( g, h,Y_{+}^{r}\mathbf{\phi})=\Bigg(\frac{\xi(2)L(1/2,\Sym^{2}\tau\times \pi)}{L(1,\pi,\Ad)L(1,\tau,\Ad)}\Bigg)^{-1}\frac{\int_{0}^{\infty}(\cosh t)^{-(2\kappa'+2)}2\pi \sinh 2tdt}{||v_{r}^{\hol}||^{2}}.
\end{equation}
By definition,
\begin{align*}
 \frac{\xi(2)L(1/2,\Sym^{2}\tau\times \pi)}{L(1,\pi,\Ad)L(1,\tau,\Ad)}&=\frac{\pi^{-1}\Gamma(1)
 \cdot
 2^{2}(2\pi)^{-2\kappa'-1}\Gamma(\kappa'+\kappa)\Gamma(\kappa'-\kappa+1)2(2\pi)^{-\kappa}\Gamma(\kappa)}{2(2\pi)^{-\kappa'-1}\Gamma(\kappa'+1)\pi^{-1}\Gamma(1)
 \cdot
 2^{2}(2\pi)^{-2\kappa}\Gamma(2\kappa)\pi^{-1}\Gamma(1)}\\
&=(2\pi)^{1-r}\frac{r!}{r+\kappa}\frac{{r+2\kappa-1 \choose r}}{{r+\kappa-1\choose r}}.
\end{align*}
Moreover, 
\begin{equation*}
    ||v_{r}^{\hol}||^{2}=\sum_{\substack{0\leq l\leq r \\ l \in 2\mathbb{N}_{0}}}|c_{l}|^{2}||v_{r-l,l}||^{2}
    =1+\sum_{\substack{2\leq l\leq r \\ l \in 2\mathbb{N}_{0}}}\prod_{\substack{0\leq j\leq l-2 \\ j \in 2\mathbb{N}_{0}}}\frac{(r-j)(r-j-1)}{(j+2)(2\kappa+j+1)}=2^{r}\frac{{r+\kappa-1\choose r}}{{r+2\kappa-1\choose r}}
\end{equation*}
The first two equalities above follow from (\cite{Xue19}, Lemma 3.3) and the last one from Lemma \ref{Monsterid1} in the appendix.
It is also well-known that,
$$
\int_{0}^{\infty}(\cosh t)^{-(2\kappa'+2)}2\pi \sinh 2tdt=4\pi(2\kappa')^{-1}.$$
Proposition \ref{archk'>k} now follows by plugging in the above explicit values into (\ref{someeqn1}).
\end{subsection}
\end{section}

\begin{section}{Computation at Nonarchimedean places}\label{nonarchplaces}
\begin{subsection}{Odd residue characteristic places}\label{oddres}
Recall the following notation. For simplicity, we suppress all subscripts $v$.
\begin{itemize}[noitemsep,topsep=0pt]
    \item $F$ is a nonarchimedaen local field with ring of integers $\mathfrak{o}$ and uniformizer $\varpi$. Let $q=\# \mathfrak{o}/(\varpi)$ and let $e$ be the integer such that $|2|^{-1}=q^{e}$.
    \item $\psi$ is a nontrivial additive character of $F$ with conductor $\delta$, i.e $\delta^{-1}\mathfrak{o}$ is the largest subgroup of $F$ on which $\psi$ is trivial. For each $a\in F^{\times}$, we denote by $\gamma_{\psi}(a)$ the Weil constant.
    \item $\tau$ is an unramified principal series representation of $\PGL_{2}(F)$ and contains a $\PGL_{2}(\mathfrak{o})$-fixed element $g$. We will assume that $g$ is normalized such that $g(1)=1$. Put 
    $$g^{(\delta)}=\tau\begin{pmatrix} 
\delta^{-1} &  \\
 & 1 
\end{pmatrix}g.$$
    \item $B(F)$ is the Borel subgrop of $\SL_{2}(F)$ consisting of upper triangular matrices of the form $n(b)m(a)$, where
    \begin{equation*}
        m(a)=\begin{pmatrix} 
a &  \\
 & a^{-1} 
\end{pmatrix},\ \ \ \
n(b)=\begin{pmatrix} 
1 & b \\
 & 1 
\end{pmatrix}; \ \ \ \ (a\in F^{\times},\ b\in F)
    \end{equation*}
Let $\widetilde{B(F)}$ be the inverse image of $B(F)$ in $\widetilde{\SL_{2}(F)}$.   
Then $\sigma=\text{Ind}_{\widetilde{B(F)}}^{\widetilde{\SL_{2}(F)}}\chi_{\overline{\psi},s}$ is the unramified principal series representation of $\widetilde{\SL_{2}(F)}$ induced from the character
\begin{equation*}
    \chi_{\overline{\psi},s}((n(b)m(a),\epsilon)=\epsilon \frac{\gamma_{\overline{\psi}}(1)}{\gamma_{\overline{\psi}}(a)}|a|^{s}.
\end{equation*}
of $\widetilde{B(F)}$ for some purely imaginary complex number $s$. Thus the underlying vector space $\widetilde{I}_{\overline{\psi}}(s)$ of $\sigma$ is the space of complex-valued functions $f$ on $\widetilde{\SL_{2}(F)}$ such that
\begin{equation*}
    f((n(b)m(a),\epsilon)g)=\epsilon\frac{\gamma_{\overline{\psi}}(1)}{\gamma_{\overline{\psi}}(a)}|a|^{s+1}f(g)\ \ \ \ \ \ (g\in \widetilde{\SL_{2}(F)}),
\end{equation*}
and $\widetilde{\SL_{2}(F)}$ acts by right translation on $\widetilde{I}_{\overline{\psi}}(s)$. It is a unitary representation with inner product
\begin{equation*}
    (f,f')=\int_{x\in F}f(w(1)n(x))\overline{f'(w(1)n(x))}dx,
\end{equation*}
where $w(1)=\begin{pmatrix} 
 & -1 \\
 1 &  
\end{pmatrix}$.
The representation $\sigma$ contains a distinguished vector, $h$, which we define below.   
    
\begin{itemize}
    \item  For any fractional ideals $\mathfrak{a}$ and $\mathfrak{b}$ of $F$, put
$$
    \Gamma[\mathfrak{a},\mathfrak{b}]=\Bigg\{\begin{pmatrix} 
 a& b \\
 c & d\end{pmatrix}
 \Big| a,d\in \mathfrak{o}, b\in \mathfrak{a},c\in \mathfrak{b}\Bigg\} .
$$
For any subgroup $\Gamma\ \text{of}\ \SL_{2}(F)$, we denote by $\widetilde{\Gamma}$ its inverse image in $\widetilde{\SL_{2}(F)}.$ Fix $\Gamma=\Gamma[\delta^{-1}\mathfrak{o},4\delta\mathfrak{o}]$. There exists a genuine character $\epsilon:\widetilde{\Gamma}\rightarrow \mathbb{C}^{\times}$ such that $\omega_{\overline{\psi}}(\gamma)\phi_{0}=\epsilon(\gamma)^{-1}\phi_{0}$, where $\gamma\in\widetilde{\Gamma}$ and $\phi_{0}$ is the characteristic function of $\mathfrak{o}$ (\cite{HI13}, Lemma 1.1). Consider the Hecke algebra $\widetilde{\mathcal{H}}=\widetilde{\mathcal{H}}(\widetilde{\Gamma}\backslash\widetilde{\SL_{2}(F)}/\widetilde{\Gamma};\epsilon)$ which is the space of compactly supported genuine functions $\phi$ on $\widetilde{\SL_{2}(F)}$ such that $\phi(\gamma_{1}X\gamma_{2})=\epsilon(\gamma_{1})\epsilon(\gamma_{2})\phi(X)$, for $\gamma_{1},\ \gamma_{2} \in \widetilde{\Gamma}$ and $X\in\widetilde{\SL_{2}(F)}$ together with a convolution product. The Hecke algebra $\widetilde{\mathcal{H}}$ acts on the space $\widetilde{I}_{\overline{\psi}}(s)$ by 
$$\sigma(\phi)f(X)=\displaystyle\int_{\widetilde{\SL_{2}(F)}}f(XX')\phi(X')dX'.$$
\item Define
$$E^{K}(X)=\begin{cases}
q^{e}\langle \phi_{0},\omega_{\overline{\psi}}(\gamma)\phi_{0}\rangle & \text{if}\ X\in \widetilde{\Gamma}[(4\delta\mathfrak{o})^{-1},4\delta\mathfrak{o}]\\
0 & \text{otherwise}.
\end{cases}
$$
Then by (\cite{HI13}, \S3  and \S6), $E^{K}$ is an idempotent in $\widetilde{\mathcal{H}}$.
\item 
By Prop. 4.6 in \cite{HI13}, there is a unique (up to scalars) element $h\in \sigma$ fixed by the idempotent $E^{K}$. Its restriction to $\widetilde{K}=\widetilde{\SL_{2}(\mathfrak{o})}$ is given by (cf. \cite{HI13}, Prop. 4.4 and p. 1982)
$$h(X)=\langle\omega_{\overline{\psi}}(Xw(2))\phi_{0},\phi_{0}\rangle,$$
where $w(2)=\begin{pmatrix} 
 & -2^{-1} \\
 2 &  
\end{pmatrix}$. Put
$$h^{(2)}=\sigma(m(2))h,$$ where $m(2)=\begin{pmatrix} 
2 &  \\
 & 2^{-1} 
\end{pmatrix}.$
\end{itemize}

\item $\omega_{\psi}$ is the Schr\"{o}dinger-Weil representation of $\widetilde{\J(F)}$.
\end{itemize}
In this subsection we prove the following:
\begin{prop}\label{nondyadicprop}
 $$\langle  g^{(\delta)}, g^{(\delta)}\rangle^{-1}\langle h^{(2)},h^{(2)}\rangle^{-1}\langle \phi,\phi\rangle^{-1}\int_{\SL_{2}(F)}\langle \tau(X) g^{(\delta)}, g^{(\delta)}\rangle\langle \overline{\sigma(X)h^{(2)},h^{(2)}\rangle\langle \omega_{\psi}(X)\phi,\phi\rangle} dX$$
 \begin{align*}
    &=\xi_{F}(2)\frac{L(1/2,\Sym^{2}\tau\times \pi)}{L(1,\pi,\Ad)L(1,\tau,\Ad)}.
\end{align*}
\end{prop}
We first explain how to reduce the proposition to the case $\delta=1$. Let $\psi$ be any non-trivial additive character of $F$ and let $\sigma$ be the representation of $\widetilde{\SL_{2}(F)}$ as above but with respect to $\psi$. To stress the dependence of $\sigma$ on $\psi$, we temporarily denote it by $\sigma_{\psi}$ and denote the idempotent $E^{K}$ of $\sigma_{\psi}$ by $E^{K}_{\psi}$. For any $a\in F^{\times}$, we write $\psi_{a}$ for the character $x\mapsto \psi(ax)$. Note that the conductor of $\psi_{\delta^{-1}}$ is one. 
Define an automorphism $r_{\delta}$ of $\widetilde{\SL_{2}(F)}$ as
$$r_{\delta}\Bigg(\begin{pmatrix}
a & b\\
c & d
\end{pmatrix},\epsilon\Bigg)=\begin{cases}
\Bigg(\begin{pmatrix}
a & \delta b\\
\delta^{-1} c & d
\end{pmatrix},\epsilon\Bigg) & c\neq 0\\
\Bigg(\begin{pmatrix}
a & \delta b\\
 & d
\end{pmatrix},(\delta,d)\epsilon\Bigg) & c= 0
\end{cases}
$$
Define a map
$$\imath_{\delta}: I_{\psi}(s)\rightarrow I_{\psi_{\delta^{-1}}}(s),\ \ \ h'\mapsto h'^{\delta},\ \ \ h'^{\delta}(X)=h'(r_{\delta^{-1}}(X)).$$
Then, 
\begin{enumerate}[label=(\roman*)]
    \item $\imath_{\delta}\circ\sigma_{\psi}(r_{\delta^{-1}}(X))=\sigma_{\psi_{\delta^{-1}}}(X)\circ \imath_{\delta}$
    \item If $h$ is $E^{K}_{\psi}$-fixed, then $h^{\delta}$ is also $E^{K}_{\psi_{\delta^{-1}}}$-fixed and moreover $h^{(2),\delta}=\pm h^{\delta,(2)}$.
    \item $\omega_{\psi}=\omega_{\psi_{\delta^{-1}}}\circ r_{\delta}$
    \item $\langle\imath_{\delta}h_{1},\imath_{\delta}h_{2}\rangle=|\delta^{-1}|\langle h_{1},h_{2}\rangle$ for any $h_{1},\ h_{2}\in \sigma_{\psi}$.
\end{enumerate}
These properties are not hard to verify and we leave the details to the reader.\\

\begin{remark}
    The above properties also hold for a dyadic place $v$, that is, for $v|2$.
\end{remark}
\textit{Proof of Proposition \ref{nondyadicprop}}
In the left-hand side of the formula in the Proposition, we make a change of variable $X=\begin{pmatrix}
a & b\\
c & d
\end{pmatrix}\mapsto \begin{pmatrix}
a & \delta^{-1}b\\
\delta c & d
\end{pmatrix}$. Then by the properties above we are reduced to the case of $\delta=1.$ In this case the Proposition follows from (\cite{Qiu14}, Lemma 4.4) once we verify that the $E^{K}$-fixed vector $h$ is also $\SL_{2}(\mathfrak{o})$-fixed. That this is indeed so can be checked from the formula $h(X)=\langle\omega_{\overline{\psi}}(Xw(2))\phi_{0},\phi_{0}\rangle$ for $X$ equal to $m(a),\ n(b),\ w$ (see \S\ref{SW rep}). We leave the details to the reader. The following fact will prove useful in the verification. When $F$ is a non-archimedean field of odd residue characteristic, $\psi$ is a non-trivial additive character of $F$ with conductor 1, and $a\in\mathfrak{o}^{\times}$, then $\gamma_{\psi}(a)=1$.\\ 

\end{subsection}

\begin{subsection}{Dyadic places}\label{dyadic}
Recall the following notation. For simplicity, we suppress all subscripts $v$.
\begin{itemize}[noitemsep,topsep=0pt]
    \item $F$ is a degree $d$ extension of $\mathbb{Q}_{2}$ with ring of integers $\mathfrak{o}$ and uniformizer $\varpi$. Let $q=\# \mathfrak{o}/(\varpi)$ and let $e$ be the integer such that $|2|^{-1}=q^{e}$.
    \item $\psi$ is a nontrivial additive character of $F$ with conductor $\delta$, i.e,  $\delta^{-1}\mathfrak{o}$ is the largest subgroup of $F$ on which $\psi$ is trivial. For each $a\in F^{\times}$, $\gamma_{\psi}(a)$ denotes the Weil constant. $\omega_{\psi}$ is the (local) Schr\"{o}dinger-Weil representation of $\widetilde{\J(F)}$ with inner product \begin{equation*}
    \langle\phi_{1},\phi_{2}\rangle=\displaystyle\int_{F}\phi_{1}(y)\overline{\phi_{2}(y)}dy\ \ \ \ \ \ (\phi_{1},\phi_{2}\in \mathcal{S}(F)).
\end{equation*}
    Our chosen vector in $\omega_{\psi}$ is $\phi=\mathds{1}_{\frac{1}{2}\mathfrak{o}_{v}}$.
    \item $\tau=\text{Ind}_{B(F)}^{\GL_{2}(F)}(|\cdot|^{s'},|\cdot|^{-s'})$ is the unramified principal series representation of $\GL_{2}(F)$ induced from the character of $B(F)$,
    \begin{equation*}
        \begin{pmatrix}
        a & b\\
        & d
        \end{pmatrix}\mapsto |ad^{-1}|^{s'}
    \end{equation*}
    for some purely imaginary complex number $s'$. It contains a $\GL_{2}(\mathfrak{o})$-fixed element $ g$ such that $ g(1)=1$. Put 
    $$ g^{(\delta)}=\tau\begin{pmatrix} 
\delta^{-1} &  \\
 & 1 
\end{pmatrix} g.$$
    \item $\sigma=\text{Ind}_{\widetilde{B(F)}}^{\widetilde{\SL_{2}(F)}}\chi_{\overline{\psi},s}$ is a principal series representation of $\widetilde{\SL_{2}(F)}$ acting on the space $\widetilde{I}_{\psi}(s)$ with definition same as the one given in \S\ref{nonarchplaces}. As before, $h$ denotes an $E^{K}$-fixed vector in $\sigma$. It is unique up to multiplication by a nonzero constant and its restriction to $\widetilde{\SL_{2}(\mathfrak{o})}$ is defined as before. Moreover, $h^{(2)}=\sigma(m(2))h$. 
    \end{itemize}
    
The goal of this section is to prove the following.

\begin{prop}\label{dyadicprop}
For $g^{(\delta)}, h^{(2)}$ and $\phi$ as above, we have
 $$\langle  g^{(\delta)}, g^{(\delta)}\rangle^{-1}\langle h^{(2)},h^{(2)}\rangle^{-1}\langle \phi,\phi\rangle^{-1}\int_{\SL_{2}(F)}\langle \tau(X) g^{(\delta)}, g^{(\delta)}\rangle\langle \overline{\sigma(X)h^{(2)},h^{(2)}\rangle\langle \omega_{\psi}(X)\phi,\phi\rangle} dX$$
 \begin{align*}
    &=2^{-d}\xi_{F}(2)\frac{L(1/2,\Sym^{2}\tau\times \pi)}{L(1,\pi,\Ad)L(1,\tau,\Ad)}.
\end{align*}
\end{prop}
As in the $v\nmid 2$ case, we may and will assume that $\delta=1$. Note that this time the $E^{K}$-fixed vector $h$ is not $\SL_{2}(\mathfrak{o})$-fixed; in fact, $\sigma$ cannot have an $\SL_{2}(\mathfrak{o})$-fixed vector since $\SL_{2}(\mathfrak{o})$ is not a subgroup of  $\widetilde{\SL_{2}(F)}$. But it is possible to rewrite the integral\\ $\displaystyle\int_{\SL_{2}(F)}\langle \tau(X) g, g\rangle\langle \overline{\sigma(X)h^{(2)},h^{(2)}\rangle\langle \omega_{\psi}(X)\phi,\phi\rangle} dX$ as 
\begin{equation*}
    \int_{\SL_{2}(F)}\langle\tau(X) g, g\rangle\overline{\langle \sigma\otimes\omega_{\psi}(X)\Phi,\Phi\rangle} dX
\end{equation*}
for some $J(\mathfrak{o})$\textit{-fixed} vector $\Phi$. We will now show how to do this. 

For $\lambda \in \frac{1}{2}\mathfrak{o}/\mathfrak{o}$, put $\phi_{\lambda}=\mathds{1}_{\lambda+\mathfrak{o}}\in \mathcal{S}(F)$. Moreover, for any $\lambda \in \frac{1}{2}\mathfrak{o}/\mathfrak{o}$, let $h_{\lambda}\in \widetilde{I}_{\overline{\psi}}(s)$ be such that 
\begin{equation*}
    h_{\lambda}(\widetilde{X})=\langle \omega_{\overline{\psi}}(\widetilde{X})\phi_{\lambda},\phi_{0}\rangle\ \ \ \ \ \ (\widetilde{X}\in \widetilde{K}).
\end{equation*}

Consider the subspace of $\mathcal{S}(F)$ spanned by the orthonormal vectors $\{\phi_{\lambda}| \lambda \in \frac{1}{2}\mathfrak{o}/\mathfrak{o}\}$. By Proposititon 3.3 in \cite{HI13}, the action of $\tilde{K}$ by $\omega_{\psi}$ on this space gives an irreducible representation of $\tilde{K}$ which we denote by $\omega_{\psi}^{K}$. It follows that the subspace of $\tilde{I}_{\psi}(s)$ spanned by $\{h_{\lambda}|\lambda\in \frac{1}{2}\mathfrak{o}/\mathfrak{o}\} $
is also invariant under the action of $\tilde{K}$ by $\sigma$. We denote this representation by $\sigma^{K}$. Noting that $|\frac{1}{2}\mathfrak{o}/\mathfrak{o}|=2^{d}$, it follows from Schur's orthogonality that $\{2^{d/2}h_{\lambda}| \lambda \in \frac{1}{2}\mathfrak{o}/\mathfrak{o}\}$ is an orthonormal basis of $\sigma^{K}$. 
Let $\imath:\omega^{K}_{\psi}\longrightarrow \sigma^{K}$ be the homomorphism with $\imath(\phi_{\lambda})=2^{d/2}h_{\lambda}$ for all $\lambda\in \frac{1}{2}\mathfrak{o}/\mathfrak{o}$. Then it is clear that $\imath$ is an isomorphism of representations of $\tilde{K}$. Moreover, it is an isometry.

\begin{lemma}\label{dyadiclemma}
\begin{itemize}
    \item [1.] We have,
    $$2^{d/2}\gamma_{\psi}(1)h^{(2)}=\sum_{\lambda\in \frac{1}{2}\mathfrak{o}/\mathfrak{o}}h_{\lambda}.$$
    \item [2.] We have,
    $$\int_{K}\sigma(k)h^{(2)}\otimes \omega_{\psi}(k)\phi dk=2^{-d/2}\gamma_{\psi}(1)^{-1}\sum_{\lambda\in \frac{1}{2}\mathfrak{o}/\mathfrak{o}}h_{\lambda}\otimes \phi_{\lambda}.$$
    \item [3.] The element
    $$\sum_{\lambda\in \frac{1}{2}\mathfrak{o}/\mathfrak{o}}h_{\lambda}\otimes \phi_{\lambda}\in \sigma\otimes \omega_{\psi}$$
is $\J(\mathfrak{o})$-fixed.
\end{itemize}
\end{lemma}

\begin{proof}
The first assertion follows from the following observations, $$\omega_{\psi}(w(1))\phi_{0}=2^{-d/2}\gamma_{\psi}(1)^{-1}\sum_{\lambda\in \frac{1}{2}\mathfrak{o}/\mathfrak{o}}\phi_{\lambda}\ \ \ \text{and}\ \ \ \sigma(w(1))h_{0}=h^{(2)}.$$
For the second assertion, we note that
$$\phi=\sum_{\lambda}\phi_{\lambda}.$$
Moreover, the representations $\sigma^{K}$ and $\overline{\omega^{K}_{\psi}}$ are dual to each other and $\{2^{d/2}h_{\lambda}\}$ and $\{\phi_{\lambda}\}$ form a dual basis. Together with Schur's orthogonality relation, we find that
\[
\bigg\langle \int_{K}\sigma(k)h^{(2)}\otimes \omega_{\psi}(k)\phi dk,h_{\lambda}\otimes \phi_{\mu}\bigg\rangle=
\begin{cases}
2^{-d} & \text{if}\  \lambda=\mu\\
0 & \text{if}\ \lambda\neq \mu.  
\end{cases}
\]
The same relation holds for the right hand side of the assertion as well. This proves (2).\\
By (2), $\sum_{\lambda\in \frac{1}{2}\mathfrak{o}/\mathfrak{o}}h_{\lambda}\otimes \phi_{\lambda}$ is $K$-invariant. It is $H(\mathfrak{o})$-invariant by the the formulae (\ref{weilrep}) of the Schr\"{o}dinger-Weil representation. This completes the proof of the lemma.
\end{proof}

The representation $\sigma\otimes \omega_{\psi}$ of the Jacobi group is isomorphic to the principal series representation (see \cite{BS98}, Thm. 5.4.2)
\begin{equation*}
    \text{Ind}_{B_{\J}(F)}^{\J(F)}\psi|\cdot|^{s}=\{f:\J(F)\longrightarrow \mathbb{C}\ |\  f(m(a)n(b)(0,\mu,\xi)g)=\psi(\xi)|a|^{s+3/2}f(g),\ g\in\J(F)\},
\end{equation*}
where the isomorphism is given by 
\begin{equation}
    h'\otimes \phi\mapsto (gu\mapsto h'(g)\omega_{\psi}(gu)\phi(0));\ \ \ \ h'\in\sigma,\phi\in \mathcal{S}(F),g\in \SL_{2}(F),u\in\He(F).
\end{equation}
We denote this principal series representation by $\rho$ and the underlying space of functions by $I(s,\psi)$. For a suitable choice of inner product on $\rho$, $\rho\simeq \sigma\otimes \omega_{\psi}$ is an isometry. We will often identify $\rho$ with $\sigma\otimes \omega_{\psi}$. \\
Put
$$\Phi=\sum_{\lambda\in \frac{1}{2}\mathfrak{o}/\mathfrak{o}}h_{\lambda}\otimes \phi_{\lambda}\in I(s,\psi).$$
Then $\Phi$ is $\J(\mathfrak{o})$-fixed by the lemma above. We have,

\begin{lemma}\label{rewritten}
\begin{equation*}
    \displaystyle\int_{\SL_{2}(F)}\langle \tau(X) g, g\rangle\langle \overline{\sigma(X)h^{(2)},h^{(2)}\rangle\langle \omega_{\psi}(X)\phi,\phi\rangle} dX=2^{-d}\int_{\SL_{2}(F)}\langle\tau(X) g, g\rangle\overline{\langle \rho(X)\Phi,\Phi\rangle} dX.
\end{equation*}
\end{lemma}

\begin{proof}
Let $X=k_{1}\diag[\varpi^{t},\varpi^{-t}]k_{2}$ where $k_{1},k_{2}\in K$ and $t$ is a nonnegative integer be the Cartan decomposition of $X\in\SL_{2}(F)$. We will denote $\diag[\varpi^{t},\varpi^{-t}]$ by $t$. By Lemma \ref{dyadiclemma}, we have 
\begin{align*}
    &\ \ \ \  \int_{\SL_{2}(F)}\langle \tau(X) g, g\rangle\langle \overline{\sigma(X)h^{(2)},h^{(2)}\rangle\langle \omega_{\psi}(X)\phi,\phi\rangle} dX\\
    &=\int_{\SL_{2}(F)}\langle \tau(X) g, g\rangle\langle \overline{\rho(X)(h^{(2)}\otimes\phi),h^{(2)}\otimes \phi\rangle} dX\\
    &=\int_{t}\langle \tau(t) g, g\rangle\Bigg(\int_{K}\int_{K}\overline{\langle\rho(tk_{2})(h^{(2)}\otimes\phi),\rho(k_{1}^{-1})(h^{(2)}\otimes \phi)\rangle}dk_{1}dk_{2}\Bigg)dt,\ \ \ \ \ \ \ \ (\text{$g$ is $K$-fixed})  \\
    &=2^{-d/2}\int_{t}\langle \tau(t) g, g\rangle\Bigg(\int_{K}\overline{\langle\rho(tk_{2})(h^{(2)}\otimes\phi),\Phi\rangle}dk_{2}\Bigg)dt,\ \ \ \ \ \ \ \ \ \ \ \ \ \ \ \ \ \ \ \ \ \ \ \ \ \ \  (\text{Lemma \ref{dyadiclemma} (ii)})\\
    &=2^{-d}\int_{t}\langle\tau(t) g, g\rangle\overline{\langle \rho(t)\Phi,\Phi\rangle} dt \ \ \ \ \ \ \ \ \ \ \ \ \ \ \\
    &=2^{-d}\int_{\SL_{2}(F)}\langle\tau(X) g, g\rangle\overline{\langle \rho(X)\Phi,\Phi\rangle} dX \ \ \ \ \ \ \ \ \ \ \ \ \ \ \ \ \ \ \ \ \ \ \  (\text{$\Vol(K)=1$ and $\Phi$ is also $K$-fixed})
    \end{align*}
\end{proof}
We will now compute
\begin{equation}\label{dyadicperiod}
   \int_{\SL_{2}(F)}\langle\tau(X) g, g\rangle\overline{\langle \rho(X)\Phi,\Phi\rangle} dX. 
\end{equation}
The computation is similar to that done in \S6 of \cite{Xue18}. Let 
$$\eta=(1,0,0)\begin{pmatrix}
&-1\\
1 &
\end{pmatrix}\in \J(F).$$
 Let
\begin{equation*}
    T_{s',s,\psi}=\int_{\SL_{2}(F)} g(X)\overline{\Phi(\eta X)}dX.
\end{equation*}
Then by (\cite{Xue17}, \S4), we have 
\begin{equation}
    (\ref{dyadicperiod})=\frac{\xi_{F}(2)}{\xi_{F}(1)^{2}}T_{s',s,\psi}T_{-s',-s,\psi^{-1}}
\end{equation}
Therefore, it is enough to calculate $T_{s',s,\psi}$.
\begin{lemma}
The support of $\Phi$ is contained in $B_{\J}(F)\J(\mathfrak{o})$.
\end{lemma}
\vspace{-1.3em}
\begin{proof}
Substitute $F$ for $\mathbb{Q}_{2}$ and $\mathfrak{o}$ for $\mathbb{Z}_{2}$ in the proof of (\cite{Xue18}, Lemma 6.4).
\end{proof}
Let $X=m(a)n(x)k$ be the usual Iwasawa decomposition for $X\in \SL_{2}(F)$. Then $dX=d^{\times}adxdk$, where $d^{\times}a$ is the measure on $F^{\times}$ such that $vol\ \mathfrak{o}^{\times}=1$. Then $ g(X)=|a|^{s'+1}$. Moreover,
\begin{equation*}
    \Phi(\eta m(a)n(x))=\begin{cases}
    |a|^{-s-\frac{3}{2}}\mathds{1}_{\mathfrak{o}}(a^{-1})& |x|\leq 1\\
    |ax|^{-s-\frac{3}{2}}\mathds{1}_{\mathfrak{o}}(a^{-1}x^{-1})\psi(a^{-2}x^{-1})& |x|> 1.
    \end{cases}
\end{equation*}
The case $|x|\leq 1$ is straightforward. 
The case $|x|>1$ relies on the following decomposition of matrices,
\begin{equation*}
    \eta m(a)n(x)=(1,0,0)m(a^{-1})w(1)n(x)=(1,0,0)m(a^{-1})m(x^{-1})n(x)\begin{pmatrix}
1&\\
x^{-1} &1 
\end{pmatrix}.
\end{equation*}
Therefore,
$$T_{s',s,\psi}=\int_{|x|\leq 1}\int_{|a|\geq 1}|a|^{-\frac{1}{2}-s+s'}d^{\times}adx+\int_{|x|>1}\int_{|ax|\geq 1}|a|^{-\frac{1}{2}-s+s'}|x|^{-\frac{3}{2}-s}\psi(a^{-2}x^{-1})dxd^{\times}a.$$
We have
$$\int_{|x|\leq 1}\int_{|a|\geq 1}|a|^{-\frac{1}{2}-s+s'}d^{\times}adx=\sum_{n=0}^{n=\infty}q^{n(-\frac{1}{2}-s+s')}=\frac{1}{1-q^{-\frac{1}{2}-s+s'}},$$
and 
$$\int_{|x|>1}\int_{|ax|\geq 1}|a|^{-\frac{1}{2}-s+s'}|x|^{-\frac{3}{2}-s}\psi(a^{-2}x^{-1})dxd^{\times}a$$
$$=\sum_{m=1}^{\infty}\sum_{n=-m}^{\infty}q^{n(-\frac{1}{2}-s+s')+m(-\frac{3}{2}-s)}\int_{|x|=q^{m}}\int_{|a|=q^{n}}\psi(a^{-2}x^{-1})dxd^{\times}a.$$
Since
$$\int_{x\in\mathfrak{o}^{\times}}\psi(\xi x)dx=\begin{cases}
1-q^{-1}& \xi \in \mathfrak{o}\\
-q^{-1} & \text{ord}(\xi)=-1\\
0 & \text{otherwise},
\end{cases}$$
we see that
$$\int_{|x|=q^{m}}\int_{|a|=q^{n}}\psi(a^{-2}x^{-1})dxd^{\times}a=\begin{cases}
q^{m-1}(q-1)& 2n+m\geq 0\\
-q^{m-1} & 2n+m=-1\\
0 & otherwise.
\end{cases}
$$
Therefore,
$$ \int_{|x|>1}\int_{|ax|\geq 1}|a|^{-\frac{1}{2}-s+s'}|x|^{-\frac{3}{2}-s}\psi(a^{-2}x^{-1})dxd^{\times}a  $$

\begin{align*}
    &=
    q^{-1}(q-1)\sum_{m\geq1}\sum_{n\geq -\ceil{\frac{m}{2}}}q^{n(-\frac{1}{2}-s+s')+m(-\frac{1}{2}-s)}-q^{-1}\sum_{{m\geq1},m\ \text{odd}}q^{\frac{m+1}{2}(\frac{1}{2}+s-s')+m(-\frac{1}{2}-s)}\\
&=\frac{-1+q-q^{\frac{1}{2}+s}+q^{1+s'}}{q^{\frac{3}{2}+s+s'}(1-q^{-\frac{1}{2}-s-s'})(1-q^{-\frac{1}{2}-s+s'})}.
\end{align*}

It follows that 
\begin{align*}
    T_{s',s,\psi}&= \frac{1}{1-q^{-\frac{1}{2}-s+s'}}+\frac{-1+q-q^{\frac{1}{2}+s}+q^{1+s'}}{q^{\frac{3}{2}+s+s'}(1-q^{-\frac{1}{2}-s-s'})(1-q^{-\frac{1}{2}-s+s'})}\\
    &=\frac{(1-q^{-1-s'})(1+q^{-\frac{1}{2}-s})}{(1-q^{-\frac{1}{2}-s+s'})(1-q^{-\frac{1}{2}-s-s'})}.
\end{align*}
We thus end up with
\begin{align*}
    (\ref{dyadicperiod})&=\frac{\xi_{F}(2)}{\xi_{F}(1)^{2}}\frac{(1-q^{-1-s'})(1+q^{-\frac{1}{2}-s})(1-q^{-1+s'})(1+q^{-\frac{1}{2}+s})}{(1-q^{-\frac{1}{2}-s+s'})(1-q^{-\frac{1}{2}-s-s'})(1-q^{-\frac{1}{2}+s-s'})(1-q^{-\frac{1}{2}+s+s'})}\\
    &=\xi_{F}(2)\frac{L(1/2,\Sym^{2}\tau\times \pi)}{L(1,\pi,\Ad)L(1,\tau,\Ad)}.
\end{align*}
\textit{Proof of Proposition \ref{dyadicprop}}. 
The above equation together with lemma \ref{rewritten} gives,
\begin{equation*}
    \int_{\SL_{2}(F)}\langle \tau(X) g, g\rangle\langle \overline{\sigma(X)h^{(2)},h^{(2)}\rangle\langle \omega_{\psi}(X)\phi,\phi\rangle} dX=2^{-d}\xi_{F}(2)\frac{L(1/2,\Sym^{2}\tau\times \pi)}{L(1,\pi,\Ad)L(1,\tau,\Ad)}.
\end{equation*}
We also have $\langle g, g\rangle=1$, $\langle\phi,\phi\rangle=2^{d}$ and by Lemma \ref{dyadiclemma}, $\langle h^{(2)},h^{(2)}\rangle=2^{-d}$. \qed
\end{subsection}

\end{section}

\appendix

\begin{section}{Some combinatorial identities}

 \begin{lemma}\label{gammaid}
  Let $N$ be a nonnegative integer. We have the following identity.
  \begin{equation}
      \sum_{i=0}^{n}(-1)^{n}\binom{N}{i}\frac{\Gamma(z+i)}{\Gamma(w+i)}=\frac{\Gamma(z)}{\Gamma(w-z)}\cdot\frac{\Gamma(w-z+N)}{\Gamma(w+N)},
  \end{equation}
  for every $z,w\in\mathbb{C}$.
 \end{lemma}
\begin{proof}
Lemma 2.1 in \cite{Ik98}.
\end{proof}
\begin{lemma}\label{Monsterid1}
For every $r,\kappa\in \mathbb{N}$, we have
 \begin{equation}\label{monsterid1}
     \sum_{0\leq l\leq \floor{\frac{r}{2}}}\prod_{0\leq j\leq l-1}\frac{(r-2j)(r-2j-1)}{(2j+2)(2\kappa+2j+1)}=2^{r}\frac{{r+\kappa-1\choose r}}{{r+2\kappa-1\choose r}},  
 \end{equation}
 where we take the empty product corresponding to $i=0$ to have value one.
\end{lemma}
\begin{proof}
We only sketch the proof here. Through some tedious but elementary calculations, one can verify that
\begin{equation*}
    (\ref{monsterid1})=\frac{\Gamma(\kappa+\frac{1}{2})}{\Gamma(\frac{\epsilon-r}{2})}\sum_{0\leq l\leq \floor{\frac{r}{2}}}(-1)^{l}\binom{\floor{\frac{r}{2}}}{l}\frac{\Gamma(l+\frac{\epsilon-r}{2})}{\Gamma(l+\kappa+\frac{1}{2})},
\end{equation*}
 where $\epsilon\equiv 1+r (\text{mod} 2)$.  Now use Lemma \ref{gammaid}.
\end{proof}

\begin{lemma}\label{Monsterid2}
For every $r,\kappa'\in \mathbb{N}$, we have
 \begin{equation} \label{monsterid2}
    \sum_{0\leq i\leq \floor{\frac{r}{2}}}\prod_{0\leq l\leq i-1} \frac{(r-2l)^{2}(r-2l-1)^{2}}{(2i-2l)(2i+2\kappa'-2l)(r-2i+2l+1)(r-2i+2l+2)} =2^{-r}\frac{{2\kappa'+2r\choose r}}{{\kappa'+r \choose r}},
 \end{equation}
 where we take the empty product corresponding to $i=0$ to have value one.
\end{lemma}
\begin{proof}
We only sketch the proof here. Through some tedious but elementary calculations, one can verify that
\begin{equation*}
    (\ref{monsterid2})=\frac{\Gamma(\kappa'+1)}{\Gamma(\frac{\epsilon-r}{2})}\sum_{0\leq l\leq \floor{\frac{r}{2}}}(-1)^{i}\binom{\floor{\frac{r}{2}}}{i}\frac{\Gamma(i+\frac{\epsilon-r}{2})}{\Gamma(i+\kappa'+1)},
\end{equation*}
where $\epsilon\equiv 1+r (\text{mod} 2)$. Now use Lemma \ref{gammaid}.
\end{proof}
\end{section}

\end{document}